\documentclass{article}
\usepackage{amsmath}
\usepackage{amsfonts}
\usepackage{bm}
\usepackage{graphicx}
\usepackage[a4paper, total={6in, 8in}]{geometry}
\usepackage[colorinlistoftodos]{todonotes}
\usepackage{hyperref}
\hypersetup{colorlinks=true, allcolors=blue}
\usepackage{float}
\usepackage{cases}
\usepackage{amssymb}
\usepackage{amsbsy}
\usepackage{bm}
\usepackage{siunitx}
\usepackage[toc,page]{appendix}
\usepackage{psfrag}
\usepackage{pdfpages}
\usepackage{mathrsfs}
\usepackage{wrapfig}
\numberwithin{equation}{section}

\newcommand{\me}{\mathrm{e}}
\newcommand{\mi}{\mathrm{i}}
\definecolor{ao(english)}{rgb}{0.0, 0.5, 0.0}
\usepackage{caption}
\usepackage{subcaption}
\usepackage{upgreek}
\usepackage[utf8]{inputenc}

\newcommand{\bbR}{\mathbb{R}}
\renewcommand{\a}{\mathbf{a}}
\newcommand{\x}{\mathbf{x}}
\newcommand{\dd}{\mathbf{d}}
\newcommand{\bdelta}{\boldsymbol{\delta}}
\newcommand{\inc}{\text{inc}}
\newcommand{\obs}{\text{obs}}
\newcommand{\spline}{\text{spline}}
\newcommand{\f}{\mathbf{f}}
\newcommand{\xh}{\widehat{\mathbf{x}}}

\begin{document}

\title{A surrogate-Bayesian algorithm for scatterer shape identification from phaseless data}

\author{Erik~Garc\'ia~Neefjes${}^{1}\footnote{E-mail: \texttt{erik.garcia@mq.edu.au}}$,
  Stuart C.~Hawkins${}^{1}\footnote{Author for correspondence. E-mail: \texttt{stuart.hawkins@mq.edu.au}}$,
Mahadevan Ganesh${}^{2}\footnote{E-mail: \texttt{mganesh@mines.edu}}$
\\[5pt]
{\footnotesize
${}^{1}$ School of Mathematical and Physical Sciences, Macquarie University, Sydney NSW 2109, Australia}
\\
{\footnotesize
${}^{2}$Department of Applied Mathematics and Statistics, Colorado School of Mines, Golden, CO 80401, USA}}

\date{\today}

\maketitle

\begin{abstract}
  This work addresses the reconstruction of a scatterer's shape
  from phaseless far field-intensity data
arising from multiple incident waves interacting
with the scatterer.
We 
  formulate the reconstruction as a statistical inverse
  scattering problem and 
  adopt a Bayesian inference framework, which
  can readily be used to compute statistical moments 
  for quantification of uncertainties in the shape
  reconstruction
  that arise from noise in the data due to measurement constraints.
  The shape of the scatterer is represented by a spline-based prior, with
  Bayesian parameters defined at the spline’s knots. To efficiently evaluate
  the Bayesian likelihood across thousands of sampling points, we develop the
  intensity  property inspired neural network (IPINN) surrogate. This surrogate
  incorporates the Helmholtz equation in the unbounded domain, 
  exterior to each sampled scatterer, along with the radiation condition at
  infinity, enabling fast and accurate simulation of the acoustic far-field
  intensity.
  Importantly, the IPINN surrogate is trained independently of the observed data
  and requires only a single incident wave for training. We demonstrate that
  this surrogate approach yields a speed-up of several orders of magnitude.
  The resulting IPINN-Bayesian
  framework offers an efficient solution for  shape reconstruction in unbounded
  domains with multiple incident wave boundary conditions, while exactly
  enforcing the radiation condition. Numerical experiments confirm the
  efficiency and effectiveness of the proposed algorithm.
\end{abstract}


\section{Introduction}\label{section:Intro}
We consider the problem of determining the shape of an object using acoustic
intensity data generated by multiple incident plane waves at a fixed frequency. This  phaseless, that is, modulus-only, far field-data problem is a class of
frequency-domain inverse scattering (FDIS) problem.
FDIS problems have numerous practical applications~\cite{colton:inverse}.
The shape-identification FDIS problem has been extensively studied
in the deterministic setting, along with the associated forward model governed by the Helmholtz partial differential equation (PDE) exterior to the scatterer
and the Sommerfeld radiation condition (SRC) at infinity~\cite{colton:inverse}. 
As with many PDE-based inverse problems, the deterministic problem is ill-posed, and regularization techniques are necessary to obtain stable solutions.

In this work, we focus on the two-dimensional sound-soft scatterer problem and address it as a statistical inference formulation of FDIS. The Bayesian framework provides a principled approach to solving PDE-based inverse problems that
enables uncertainty quantification (UQ)~\cite{stuart:inverse, dashtistuart:inverse}.
Crucially, uncertainty due to
noise, data insufficiency and geometric modeling uncertainties 
are evident in the model output and can be incorporated into decision making.
However, the high computational cost of Bayesian methods—due to repeated PDE solves in high-dimensional parameter spaces for likelihood evaluations—has historically limited their use~\cite{stuart:inverse,dashtistuart:inverse}. Applications of Bayesian inference to the shape FDIS problem with the exterior Helmholtz PDE remain scarce. Even at the forward level, rigorous analysis of shape uncertainties for the Helmholtz PDE has only recently been developed, relying on Karhunen–Lo\'eve (KL) expansions with global Fourier 
bases, see~\cite{kuijpers2024wavenumber} and references therein.

The literature on shape reconstruction from phaseless data is limited. A central difficulty arises from the translational invariance of far-field magnitudes, which complicates the establishment of uniqueness in the FDIS problem~\cite{inv_acou_2019_pap}. Indeed, no general uniqueness results are known for arbitrary scatterers, even when invariance is addressed through techniques such as incident-field superposition, and even when restricted to the two-dimensional Helmholtz PDE with a sound-soft boundary~\cite{inv_acou_2019_pap}. Numerical evidence
for the deterministic model in~\cite{inv_acou_2019_pap} and references therein
suggests that reconstructions from phaseless data are feasible.
In particular, for the two-dimensional model, sampling-based methods have been demonstrated in~\cite{inv_acou_2019_pap}, while iterative approaches have been explored in~\cite{inv_acou_2016_pap}. 
We are not aware of any reconstructions from phaseless data using
the Bayesian framework.

Early Bayesian studies of FDIS have primarily relied on KL representations. For example,~\cite{bui2014analysis} analyzed the two-dimensional inverse FDIS problem with phased data using a standard Bayesian framework and KL priors, while~\cite{yang2020bayesian, ganesh2020efficient} extended this approach to phaseless acoustic intensity data. In contrast, we propose a local spline-based prior representation of the scatterer, where the spline knot values serve as Bayesian parameters. This yields an $N_\text{spline}$-dimensional Bayesian model informed by 
$N_\text{obs} \times N_\text{inc}$ far-field measurements, with $N_\text{spline}$, $N_\text{inc}$, and $N_\text{obs}$ denoting the number of spline knots, incident waves, and observation directions, respectively.
  Numerical results in the literature~\cite{palafox2017point, neefjes2025far}
  demonstrate that it is sometimes important to use 
data from large numbers of incident waves $N_\inc$. 
The computational bottleneck in such high-dimensional models lies in repeatedly evaluating the Bayesian likelihood function and associated forward model
over many experimental configurations.

To overcome this challenge, we develop a surrogate model based on a {\bf neural network (NN)} that learns the far-field intensity response from the spline-based parameterization. Specifically, we introduce the {\bf Intensity Property Inspired NN (IPINN)}, which approximates the far-field map from the $N_\text{spline}$-dimensional parameter space to the reduced $N_\text{obs}$-dimensional observation space. By exploiting the rotational invariance of intensity, IPINN requires only a single reference incident wave and applies cyclic permutation operators to recover responses for other appropriate directions.
Thus our approach is particularly efficient when the number 
of incident waves $N_\inc$ is large.
Importantly, this surrogate is {\bf data-independent} and trained offline, enabling fast online Bayesian likelihood evaluations. This leads to our proposed IPINN-Bayesian algorithm, an offline/online framework that combines deep learning of the far-field map with efficient Bayesian inference for shape reconstruction.

In our numerical experiments we explore the computed Bayesian posterior
distribution using Markov Chain Monte Carlo (MCMC) sampling.
Our IPINN surrogate can be used in any of the standard MCMC algorithms
for forward model evaluation.
Since the NN is readily differentiable with respect to its
inputs, the IPINN can provide the gradient information required
for methods such as Hamiltonian Monte Carlo (HMC).
We present results obtained using the Gibbs sampling algorithm.

The rest of this article is organized as follows.
In Section~\ref{section:inverse-shape problem}
we describe our inverse-shape FDIS problem, including identifying our prior
space.
In Section~\ref{section: forward computation}
we state the governing PDEs and describe the numerical method that we
use for generating training data for our surrogate model.
In Section~\ref{sec:nn} we described our NN-based surrogate model.
In Section~\ref{sec: recentering} we describe how we sample the posterior
distribution and manage the phaseless data using a post-processing
scheme to align the shapes.
Finally, in Section~\ref{section:Numerical Results} we demonstrate the
efficiency and accuracy of our method for a range of test scatterers.

\section{The inverse-shape FDIS problem}\label{section:inverse-shape problem}
In this section, we formulate the
inverse-shape FDIS problem
of recovering the shape of a closed star-shaped sound-soft scatterer in
two dimensions
as a probabilistic inference problem.
The inference is based on noisy, phaseless far-field intensity data generated by multiple incident plane waves interacting with the scatterer in free-space. 
We begin by introducing the mathematical notation and modeling assumptions related to the scatterer and the associated input and output data.

Let $D(\boldsymbol{\xi})$ denote the unknown sound-soft scatterer, with its boundary parametrized by a vector $\boldsymbol{\xi}$ defined on  a probability space $(\Xi,\mathcal{F},P)$, where 
$\Xi$ is the sample space, $\mathcal{F}$ the  event space  and $P$ the probability measure. Our goal is to construct a posterior probability distribution for $\boldsymbol{\xi}$, conditioned on the phaseless data obtained from
measurements of the far fields induced by 
the scatterer interacting with several incident plane waves
with fixed wavenumber $k$ and distinct directions $\widehat{\mathbf{d}}$.
This posterior distribution yields
statistical information about
the scatterer boundary $\partial D(\boldsymbol{\xi})$, which we
extract using MCMC sampling.
The mean shape, approximated from the MCMC samples, provides an approximation
to the scatterer's shape.

To make the problem well-posed for a practical Bayesian setting, we assume that the scatterer boundary $\partial D(\boldsymbol{\xi})$ lies entirely within an annular region centered at the origin, bounded by an inner radius $R_0$ and an outer radius $R$. This assumption provides an isotropic prior representation, which avoids introducing directional biases into the reconstruction process.
A key ingredient in the formulation is the mapping $\boldsymbol{\xi} \mapsto \partial D(\boldsymbol{\xi})$,  which  defines the shape parametrization of the scatterer.

The star-shaped property of the two-dimensional
scatterer guarantees that the boundary can be represented in polar coordinates as
\begin{equation}
    \bm{\chi}(\theta) = r(\theta;\boldsymbol{\xi}) \, \widehat{\mathbf{x}}(\theta), \qquad \theta \in [0, 2\pi),
    \label{scatterer polar representation}
\end{equation}
where $\bm{\chi}(\theta) \in \mathbb{R}^2$ denotes a boundary point on the scatterer, and $\widehat{\mathbf{x}}(\theta)= (\cos \theta,\sin \theta)$ is the unit radial vector on the 
unit circle  $\partial B \subseteq \mathbb{R}^2$.  Closure of the boundary requires
$\bm{\chi}(0)=\bm{\chi}(2 \pi)$.
The positive radial function $r(\theta;\boldsymbol{\xi})$  is  represented in logarithmic form as 
\begin{equation}
    {r}(\theta;\boldsymbol{\xi}) =  e^{s(\theta;\boldsymbol{\xi})}.
    \label{eq:logradius}
\end{equation}
We define the prior space for our Bayesian formulation as the set of admissible shapes for which the log-radius $s(\theta;\boldsymbol{\xi})$
has distribution $\mathcal{U}(\log R_0,\log R)$ for $\theta \in [0,2\pi)$, where throughout this article $\log$ denotes the natural logarithm.
This ensures that the boundary remains within the prescribed annular region.


We recall that the far-field, denoted throughout this article as $u^\infty$, 
is a function defined on $\partial B$ the set of all possible observed unit vector directions  $\widehat{\mathbf{x}}$. In the FDIS model with fixed wavenumber $k$, 
the far field of a chosen scatterer
crucially depends on two factors: the scatterer boundary $\partial D$ and the incident plane wave direction $\widehat{\mathbf{d}}$.  Accordingly, to describe the
phaseless intensity
of the far-field observed at $\widehat{\mathbf{x}}$,   it is convenient to introduce the notation
\begin{equation}
\mathcal{F}({\partial D};\widehat{\mathbf{d}})(\widehat{\mathbf{x}})
     = | u^\infty(\widehat{\mathbf{x}};{\partial D};\widehat{\mathbf{d}})|^2,  \qquad \widehat{\mathbf{x}} \in \partial B,  \label{Forward operator intro}
\end{equation}
for the  nonlinear far-field \textit{forward model} operator.

In practice, intensity data is collected for only a finite number $N_\text{inc}$ of incident directions
$\dd_1,\dots,\dd_{N_\inc}$
and a finite number  $N_\text{obs}$  of observed directions 
$\x_1,\dots,\x_{N_\obs}$.
Thus the real-valued noisy intensity data
is $\bdelta=(\bdelta_1,\dots,\bdelta_{N_\inc})$
where $\bdelta_1,\dots,\bdelta_{N_\inc}$
are
vectors of length $N_\text{obs}$ corresponding to the $N_\inc$ incident waves.
The corresponding forward model output is $\f = (\f_1,\dots,\f_{N_\inc})$
where $\f_j = (f_{1,j},\dots,f_{N_\obs,j})$
and $f_{i,j} =
\mathcal{F}({\partial D};\widehat{\mathbf{d}_j})(\widehat{\mathbf{x}}_i)$, with $i = 1,\dots,N_\obs$ and $j = 1,\dots,N_\inc$. 
We denote the noise in the data by $\boldsymbol{\eta} \in \mathbb{R}^{N_\text{obs} \times N_{\text{inc}}}$ so that
 $\boldsymbol{\delta} =  \mathbf{f} + \boldsymbol{\eta}$.

Our stochastic model is built on the 
assumption that the noise is random and independent at each observation
direction. Then
$\boldsymbol{\eta}$ is modeled as a sample of independent and identically distributed (i.i.d.)~normal random variables with zero mean and variance $\sigma^2$, that is, $\eta_{i,j} \sim \mathcal{N}(0,\sigma^2)$ for $i=1,\dots,N_{\text{obs}}$
and $j=1,\dots,N_{\text{inc}}$. 
Since, in general, the standard deviation $\sigma$ is not known, we adopt a hierarchical model in which $\sigma$ is itself modeled as a random variable whose value is to be inferred from the data. 


Using the stochastic parametrization $\boldsymbol{\xi}$
introduced earlier to represent the scatterer, the stochastic inverse-shape FDIS problem is: conditioned on the noisy data 
$\boldsymbol{\delta}$, compute a joint posterior distribution for
$\boldsymbol{\xi}$ and $\sigma$. 
We use Bayes' Theorem to obtain the  data-conditioned \textit{posterior} distribution with density
\begin{equation}
    p(\boldsymbol{\xi}, \sigma \mid \boldsymbol{\delta}) \propto p(\boldsymbol{\delta}\mid\boldsymbol{\xi}, \sigma) \, p_{\boldsymbol{\xi}}(\boldsymbol{\xi}) \,  p_\sigma({\sigma}).\label{Bayes Thm}
\end{equation}
Here $p_{\boldsymbol{\xi}}(\boldsymbol{\xi})$ and $p_\sigma({\sigma})$
are the densities of the 
data-independent prior distributions for
$\boldsymbol{\xi}$ and $\sigma$ respectively.


For fixed observational data $\boldsymbol{\delta}$, the \textit{likelihood} function $p(\boldsymbol{\delta}\mid\boldsymbol{\xi}, \sigma)$ in (\ref{Bayes Thm}) follows from the Gaussian assumption on the noise $\boldsymbol{\eta} = \boldsymbol{\delta} -\mathbf{f}(\boldsymbol{\xi})$, so that
\begin{equation}
  p(\boldsymbol{\delta}\mid\boldsymbol{\xi}, \sigma)
  = \frac{1}{(2\pi\sigma^2)^\frac{N_{\text{obs}}}{2}}
  \exp\left(
  -\frac{(\boldsymbol{\delta} - \mathbf{f})^T(\boldsymbol{\delta} - \mathbf{f})}{2\sigma^2}\right).
    \label{eq:likelihood_normal}
\end{equation}
Using notation that explicitly indicates the dependence of the
scatterer boundary $\partial D$ on $\boldsymbol{\xi}$, we write
$\mathbf{f}(\boldsymbol{\xi})
= (\f_1(\boldsymbol{\xi}),\dots,\f_{N_\inc}(\boldsymbol{\xi}))
\in \mathbb{R}^{N_\text{inc} \times N_{\text{obs}}}$ with components 
\begin{equation}\label{eq:likely-far}
f_{i,j}(\boldsymbol{\xi}) = \mathcal{F}(\partial D (\boldsymbol{\xi});\widehat{\mathbf{d}}_j )(\widehat{\mathbf{x}}_i)=|u^\infty(\widehat{\mathbf{x}}_i;\partial D (\boldsymbol{\xi}) ;\widehat{\mathbf{d}_j})|^2, \quad \text{for} \quad i=1,\dots,N_{\text{obs}},~~   j=1,\dots,N_{\text{inc}},
\end{equation}
which
we evaluate by solving the Helmholtz PDE, together with the SRC and the sound-soft boundary condition.
Where required, we solve the PDE using
a boundary integral equation reformulation that we solve to
high accuracy using a high-order Nystr\"om discretization~\cite{colton:inverse}.
High order convergence of the Nyst\"om discretization is conditional on
the boundary $\partial D(\boldsymbol{\xi})$ being at least $C^2$-smooth~\cite{colton:inverse}.
Next we consider description of the prior distribution required in~\eqref{Bayes Thm} motivated by the desired  smoothness choice.


Returning to the assumption that  the scatterer boundary is contained within an annular region centered at the origin, 
we parametrize the log-radius using the lowest degree  $C^2$-smooth splines, the smoothest cubic splines interpolating at
equally spaced knots $\vartheta_\ell$ for $\ell=1,\dots,N_\text{spline}$.
Then the 
boundary is uniquely determined by
the values of the log-radius at the knots, so that we may
parametrize the boundary as $\partial D (\boldsymbol{\xi})$  where
$\boldsymbol{\xi} = (\xi_1,\dots,\xi_{N_\text{spline}})^\top$ and
\begin{equation}\label{eq:spl-param}
\xi_\ell=  s(\vartheta_\ell;\boldsymbol{\xi}), \qquad \ell = 1,\dots,N_\text{spline}.
\end{equation}
The numerical results in Section \ref{section:Numerical Results} demonstrate that cubic splines can efficiently represent a wide variety of shapes with relatively few parameters.

Under the assumption of isotropy in our prior,
we assign identical independent uniform distributions to each spline parameter, so that $\xi_\ell \sim \mathcal{U}(\log R_0, \log R)$ for $\ell=1,\dots,N_\text{spline}$. The radius $r_\ell$ at the $\ell$th knot satisfies $r_\ell \in [R_0, R]$ as required.
Some splines ``overshoot'' the target annular region
by having some part of the curve between the knots lying outside the annulus.

\begin{figure}
    \centering
    \includegraphics[scale=0.3]{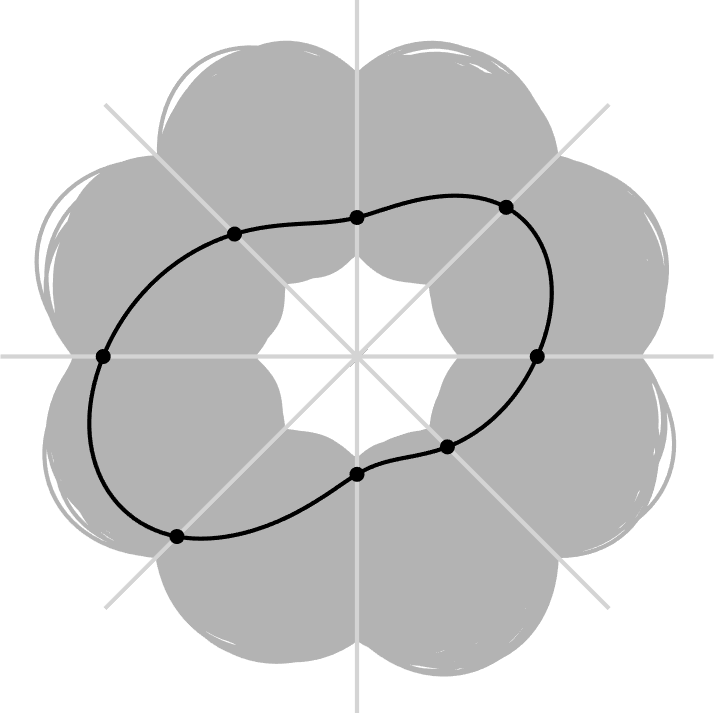}
    \caption{Visualization of $1000$ samples of the prior space (shaded grey) with $N_\text{spline}=8$ and $\xi_i \sim \mathcal{U}[-0.5,0.5]$ for $i=1,\dots,8$, so that $- \log R_0 = \log R = 0.5$. One sample is highlighted in black with its knot values marked.}
    \label{Fig:Prior Space}
\end{figure}

For the noise parameter $\sigma$, we use a log-uniform prior distribution $\log \sigma \sim \mathcal{U}[a_0^*,a^*]$
for some real parameters $a_0^*,a^*$ with $a_0^*<a^*$.
Combining these priors, the joint prior density in~(\ref{Bayes Thm}) is therefore
\begin{equation}
    p_{\boldsymbol{\xi}}(\boldsymbol{\xi}) p_\sigma({\sigma}) = \frac{1}{(R-R_0)^{N_\text{spline}}} \times \frac{1}{\sigma \log{(a^*/a_0^*)}}.\label{eqn: prior}
\end{equation}
To illustrate the boundaries described by our prior we visualize
a large number of samples in Figure~\ref{Fig:Prior Space}.

\section{Forward Model and NN training data}\label{section: forward computation}

In this section we briefly describe the BIE-based method that we use for
numerical solution of the forward model
for each fixed $\boldsymbol{\xi}$ to generate the
training data required to construct
our NN-based surrogate, which we
use for rapid evaluation of the forward model inside the MCMC sampling
process.
Training the NN requires large amounts of data, which comprises samples
of
$\boldsymbol{\xi}$ and the 
corresponding forward model data
$f_{i,j}(\boldsymbol{\xi})$ from~\eqref{eq:likely-far}  for $i=1,\dots,N_{\text{obs}}$ and $j=1,\dots,N_{\text{inc}}$.
We expand on how the training data samples are chosen in the next section.

The scattered-field $u^{\text{sc}}$ at the heart of the forward
model is induced by the plane wave 
\begin{equation}
    u^{\text{inc}}(\mathbf{x}) = \me^{\mi k \mathbf{x} \cdot \mathbf{\widehat{{d}}}},
\end{equation}
and, following details in~\cite{colton:inverse},
satisfies the exterior Helmholtz PDE,
\begin{equation}
    \left( \Delta  + k ^2 \right) u^{\text{sc}}(\mathbf{x}; \boldsymbol{\xi}) = 0, \quad \mathbf{x} \in \mathbb{R}^2 \setminus \overline{D(\boldsymbol{\xi})},
    \label{HH eqn}
\end{equation}
the sound-soft boundary condition
\begin{equation}
{u}^{\text{sc}}(\mathbf{x}; \boldsymbol{\xi}) = - u^{\text{inc}}(\mathbf{x}), \quad \text{for} \quad \mathbf{x} \in \partial D(\boldsymbol{\xi}),
\label{Dirichlet BC}    
\end{equation}
and the SRC at infinity, 
 \begin{equation}
    \lim_{|\mathbf{x}| \to \infty} \sqrt{|\mathbf{x}|} \left( \frac{\partial u^{\text{sc}}}{\partial \mathbf{x}}(\mathbf{x}; \boldsymbol{\xi}) - \mi k u^{\text{sc}}(\mathbf{x}; \boldsymbol{\xi}) \right) = 0,
    \label{2D radiation condition}
\end{equation}
uniformly with respect to the direction $\widehat{\mathbf{x}} = \mathbf{x}/ |\mathbf{x}|$. 
The associated far-field ${u}^\infty (\widehat{\mathbf{x}}; \boldsymbol{\xi})$
is related to the scattered field
${u}^{\text{sc}}(\mathbf{x}; \boldsymbol{\xi})$ by
\begin{equation}
{u}^{\text{sc}}(\mathbf{x}; \boldsymbol{\xi}) = \frac{\me^{\mi k |\mathbf{x}|}}{\sqrt{|\mathbf{x}|}} \left[{u}^\infty (\widehat{\mathbf{x}}; \boldsymbol{\xi}) + O \left( \frac{1}{|\mathbf{x}|} \right) \right], \quad \text{as} \quad |\mathbf{x}| \to \infty. 
\label{definiton u infty}
\end{equation}

For each fixed incident field,
we represent the 
scattered field using a combined field ansatz~\cite{colton:inverse},
\begin{equation}
    u^{\text{sc}}(\mathbf{x}; \boldsymbol{\xi}) = \int_{\partial D(\boldsymbol{\xi})} \left[ \frac{\partial G}{\partial \mathbf{n}(\mathbf{y})}(\mathbf{x},\mathbf{y}) - \mi k G(\mathbf{x},\mathbf{y}) \right] \phi(\mathbf{y};\boldsymbol{\xi}) \, ds(\mathbf{y}), \quad \mathbf{x} \in \mathbb{R}^2 \setminus \overline{D(\boldsymbol{\xi})},
    \label{u potentials ansatz}
\end{equation}
where the density $\phi$ is to be determined. Here 
$G(\mathbf{x}, \mathbf{y}) = \mi H_0^{(1)}(k |\mathbf{x} - \mathbf{y}|)/4$
is the free-space Green's function for the two-dimensional Helmholtz equation, $H_0^{(1)}$ is the zeroth order Hankel function of the first kind,
and $\mathbf{n}(\mathbf{y})$ represents the outward unit normal to $\partial D(\boldsymbol{\xi})$ at $\mathbf{y}$.
The unknown density $\phi$ satisfies the second kind integral
equation~\cite[Equation~(3.29)]{colton:inverse}
\begin{equation}
\phi(\mathbf{x}) + 2 \int_{\partial D(\boldsymbol{\xi})} \left[ \frac{\partial G}{\partial \mathbf{n}(\mathbf{y})}(\mathbf{x}, \mathbf{y}) - \mi k G(\mathbf{x}, \mathbf{y}) \right] \phi(\mathbf{y};\boldsymbol{\xi}) \, ds(\mathbf{y}) = -2u^{\text{inc}}(\mathbf{x}), \quad \mathbf{x} \in \partial D(\boldsymbol{\xi}),
\label{BIE final}
\end{equation}
which admits a unique solution for real-valued wavenumbers $k>0$~\cite{colton:inverse}.
After solving~\eqref{BIE final} to compute $\phi$,
we can obtain the far-field pattern explicitly via the boundary integral 
\cite[Equation~(3.110)]{colton:inverse}
\begin{equation}
    u^\infty(\widehat{\mathbf{x}};\boldsymbol{\xi}) = \frac{\me^{\mi\pi/4}}{\sqrt{8\pi k}} \int_{\partial D(\boldsymbol{\xi})} \left[ \frac{\partial e^{-\mi k \widehat{\mathbf{x}} \cdot \mathbf{y}}}{\partial \mathbf{n}(\mathbf{x}_0)} - \mi k \me^{-\mi k \widehat{\mathbf{x}} \cdot \mathbf{x}_0} \right] \phi(\mathbf{y};\boldsymbol{\xi}) \, ds(\mathbf{y}), \quad  \widehat{\mathbf{x}} \in \partial B.
    \label{far field BIE}
\end{equation}

We discretize~\eqref{BIE final}--\eqref{far field BIE} using a high-order Nystr\"om method~\cite[Section 3.6]{colton:inverse},
which uses 
the boundary mapping in \cite[Equation~(3.113)]{colton:inverse} to parametrize the boundary integrals in~\eqref{BIE final}--\eqref{far field BIE}
on $[0,2\pi)$.
We implement the Nystr\"om method using $2n+2$ evenly spaced quadrature points on $[0,2\pi]$
and construct the boundary mapping for $\partial D(\boldsymbol{\xi})$ using~\eqref{scatterer polar representation}--\eqref{eq:logradius}. 
Discretization of~\eqref{BIE final} results in a linear system of the form
\begin{equation}
\label{eq:linear-system}
    \mathbf{A}(\boldsymbol{\xi}) \boldsymbol{\phi}(\boldsymbol{\xi}) = \mathbf{b}(\boldsymbol{\xi}),
\end{equation}
where the matrix $\mathbf{A}(\boldsymbol{\xi})$ has dimension $(2n+2) \times (2n+2)$.
The Nystr\"om method exhibits high order convergence but does not realize its
potential exponential convergence in this case due to the limited smoothness of the cubic splines. Nevertheless, we can obtain satisfactory accuracy using $n \leq 100$, which allows us to solve~\eqref{eq:linear-system}
efficiently using a direct solver.

\section{Neural network based far-field surrogate: IPINN}
\label{sec:nn}

In this section we describe our surrogate for the forward model,
which we use whenever evaluations of the forward model are required
within the MCMC sampling.
We demonstrate in the numerical results section that our surrogate
model is 24,000 times faster than the Nystr\"om method that we use
to generate the training data.
Our surrogate model can be used with any of the standard algorithms used for
MCMC sampling, including
algorithms that require the
gradient of the forward model.
The gradient is readily obtained using the derivative of the NN, which is a standard functionality of NN software.

Our forward model requires far-field intensity data for several
incident wave directions. A key feature of our surrogate is
to avoid training a NN for all incident wave directions, or several
NNs---one for each incident wave directions.
Instead we take advantage of rotational symmetries to use one single
NN for all incident wave directions.

In order to exploit rotational symmetry in our model, for any $\theta \in [0,2\pi)$
  we introduce a rotation matrix 
\begin{equation}
    R_\theta = \begin{pmatrix}
\cos \theta & -\sin \theta \\
\sin \theta & \cos \theta 
\end{pmatrix},
\end{equation}
and for subsets $A \subseteq \bbR^2$ we write
$R_\theta A = \{ R_\theta \a \, : \, \a \in A \}$.
Using the rotation $R_\theta$  to rotate the
  coordinate system in which we evaluate the forward operator
  (\ref{Forward operator intro}), we have
\begin{displaymath}
    \mathcal{F}(R_\theta \partial D(\boldsymbol{\xi});R_\theta \widehat{\mathbf{d}})(R_\theta \widehat{\mathbf{x}}) = |u^\infty(R_\theta \widehat{\mathbf{x}};R_\theta \partial D(\boldsymbol{\xi});R_\theta \widehat{\mathbf{d}})|^2.
\end{displaymath}
Clearly the forward scattering problem is invariant under such a change in
coordinate system, so
that $\mathcal{F}(R_\theta \partial D(\boldsymbol{\xi});R_\theta \widehat{\mathbf{d}})(R_\theta \widehat{\mathbf{x}}) =  \mathcal{F}(\partial D(\boldsymbol{\xi});\widehat{\mathbf{d}})(\widehat{\mathbf{x}})$.

Our goal is to perform all evaluations of our forward model, which may
involve waves having several incident directions, using a
surrogate that is trained for only a single incident direction.
To that end, we
choose a principal incident direction $\widehat{\mathbf{d}}_1$ and
for a simulation with incident direction $\widehat{\mathbf{d}}$, we 
choose the angle $\theta$ so that $R_\theta \widehat{\mathbf{d}} = \widehat{\mathbf{d}}_1$. Then
\begin{equation}
  \label{eq:interim}
  \mathcal{F}(\partial D(\boldsymbol{\xi});\widehat{\mathbf{d}})(\widehat{\mathbf{x}}) =
  \mathcal{F}(R_\theta \partial D(\boldsymbol{\xi});R_\theta \widehat{\mathbf{d}})(R_\theta \widehat{\mathbf{x}}) =
  \mathcal{F}(R_\theta \partial D(\boldsymbol{\xi});\widehat{\mathbf{d}}_1)(R_\theta \widehat{\mathbf{x}}).
\end{equation}

Let $\mathbf{f}^1 = (f^1_i)_{i=1,\dots,N_\text{obs}}$ denote the forward model
data evaluated for the reference incident direction,
with  $ f^1_i = |u^\infty(\widehat{\mathbf{x}}_i;\partial D;\widehat{\mathbf{d}}_1)|^2$.
For spline knots
$\vartheta_1,\dots,\vartheta_{N_\text{spline}}$ that are evenly spaced, and particular
rotation angles $\theta = 2 \pi j/N_\text{spline}$ for some $j\in \mathbb{Z}$, we have for $\boldsymbol{\xi} = \left({\xi}_1,\dots,{\xi}_{N_\text{spline}}\right)$,
\begin{displaymath}
R_\theta \partial D(\boldsymbol{\xi}) = \partial D(\sigma_\theta \boldsymbol{\xi}),
\end{displaymath}
where $\sigma_\theta = \sigma_S^j$ and 
$\sigma_S = \left(1\,2 \dots N_{\text{spline}}\right) \in \mathbb{Z}_{N_\text{spline}}$ generates the finite cyclic group $C_{N_{\text{spline}}}$.
Similarly, for observation directions $\widehat{\mathbf{x}}_1,\dots,\widehat{\mathbf{x}}_{N_\text{obs}}$ 
that are evenly spaced, and particular rotation angles
$\theta = 2 \pi m/N_\text{obs}$ for some $m \in \mathbb{Z}$, we have 
\begin{displaymath}
  R_{\theta} \xh_i = \xh_{i'}, \qquad
  i' = i+m \, \mbox{mod} \, N_\obs.
\end{displaymath}
It follows from~\eqref{eq:interim} that, under certain conditions on $\theta$
that are clarified below,
\begin{displaymath}
  \mathcal{F}(\partial D(\boldsymbol{\xi});\widehat{\mathbf{d}})(\widehat{\mathbf{x}}) =
  \mathcal{F}(R_\theta \partial D(\boldsymbol{\xi});\widehat{\mathbf{d}}_1)(R_\theta \widehat{\mathbf{x}}) =
  \mathcal{F}(\partial D(\sigma_\theta \boldsymbol{\xi});\widehat{\mathbf{d}}_1)(\xh_{i'}),
\end{displaymath}
which can be computed using the surrogate for the incident direction
$\widehat{\mathbf{d}}_1$.
In particular, for the above choice of $(j,m)$, the forward model output is given by
\begin{equation}
\mathbf{f}^j(\boldsymbol{\xi}) = \sigma_O^{-m}  \mathbf{f}^1(\sigma_\theta \boldsymbol{\xi}), \label{output rotation}
\end{equation}
where $\sigma_O = (1\,2 \dots N_{\text{obs}}) \in \mathbb{Z}_{N_\text{obs}}$ generates of the finite cyclic group $C_{N_{\text{obs}}}$.

Now we clarify the conditions on $\theta$ under which these symmetries
can exploited. These can be stated in a general way in terms of invariance
under rotation of the incident wave directions, observation directions
and spline knots:
\begin{subequations}\label{rotation conditions}
\begin{align}
R_\theta \{\widehat{\mathbf{d}}_1,\dots,\widehat{\mathbf{d}}_{N_\text{inc}} \} & = \{\widehat{\mathbf{d}}_1,\dots,\widehat{\mathbf{d}}_{N_\text{inc}} \},\\
R_\theta \{ \widehat{\mathbf{x}}_1,\dots,\widehat{\mathbf{x}}_{N_\text{obs}}\} & = \{\widehat{\mathbf{x}}_1,\dots,\widehat{\mathbf{x}}_{N_\text{obs}}\},\\
\theta + \{ \vartheta_1,\dots,\vartheta_{N_\text{spline}}\} & = \{\vartheta_1,\dots,\vartheta_{N_\text{spline}}\},
\end{align}
\end{subequations}
where the arithmetic in the last equation is modulo $2 \pi$.
Clearly, these hold for instance when $\theta = 2 \pi l/N_\text{inc}$ for some $l \in \mathbb{Z}$ and $N_\text{inc}$ divides $N_\text{obs}= N_\text{spline}$. 

Equation~(\ref{output rotation}) implies that if we know the mapping $\boldsymbol{\xi} \rightarrow \mathbf{f}^1(\boldsymbol{\xi})$ then we can efficiently obtain the corresponding forward map for the other incident directions 
using the symmetry relation~\eqref{output rotation}.
We demonstrate below that this introduces a significant computational efficiency to the surrogate model, as it allows us to reduce the entire multi-directional scattering problem to a single reference configuration. 
The network needs to learn only a single mapping $\mathbf{f^1}: \mathbb{R}^{N_{\text{spline}}} \to \mathbb{R}^{N_{\text{obs}}}$, which reduces the model's parameter space by a factor of $N_\text{inc}-1$ and ensures consistency across all choices of incident directions.

The NN is trained over the prior space
$[\log R_0,\log R]^{N_\spline}$
established in Section~\ref{section:inverse-shape problem}.
To adequately cover this high-dimensional training space we use
$M$ data points $\boldsymbol{\xi}$ obtained as independent samples from the
prior,
analogous to 
the Monte-Carlo method.
For each sampled parameter vector $\boldsymbol{\xi}^{(m)}$, we compute the corresponding forward model data for the fixed incident wave with direction $\widehat{\mathbf{d}}_1$ using the Nyström boundary integral method described in Section \ref{section: forward computation}.
Each training datum then consists of the pair $(\boldsymbol{\xi}^{(m)},\mathbf{f}^1({\boldsymbol{\xi}}^{(m)}))$.
NN architecture details are described in Section \ref{section:Numerical Results}.

We ensure satisfactory numerical accuracy of our data by
determining a-priori a Nyström discretization parameter $n$ that ensures chosen accuracy across the entire training dataset.
For a given $\boldsymbol{\xi}$,
the relative error in the approximate solution computed
using Nystr\"om's method with discretization parameter $n$ is estimated using
\begin{equation}\label{eq: n vs n+5 rel error}
  \text{error}_{(n)}(\boldsymbol{\xi}) \approx
\frac{\| \mathbf{f}^{1}_{(n)}(\boldsymbol{\xi}) - \mathbf{f}^1_{(n+5)}(\boldsymbol{\xi}) \|_\infty}{\|\mathbf{f}^1_{(n+5)}(\boldsymbol{\xi}) \|_\infty},
\end{equation}
where $\| \cdot \|_\infty$ is the vector norm over the data at the observation
directions and
$\mathbf{f}^1_{{(n)}}$ denotes the approximation
computed using the Nystr\"om method
with discretization parameter $n$.
We estimate the error across the prior space using a Monte Carlo
approach where we compute the maximum norm over
1000 independent samples of $\boldsymbol{\xi}$
and choose $n$ sufficiently large that this estimate is below a prescribed tolerance $\epsilon_\text{tol}$.

\section{Posterior Sampling and shape alignment}
\label{sec: recentering}

Our NN-based surrogate can be used with any of the standard
algorithms used for MCMC sampling.
We demonstrate our surrogate by sampling from the posterior distribution in (\ref{Bayes Thm}) using Gibbs sampling~\cite{voss2013introduction}
which samples in each stochastic dimension sequentially.
In each dimension we sample by evaluating the posterior at $L$ discrete points,
so that for each Gibbs sample we require a total of
$L(N_\text{spline} + 1)N_\text{inc}$
solutions of the Helmholtz PDE.

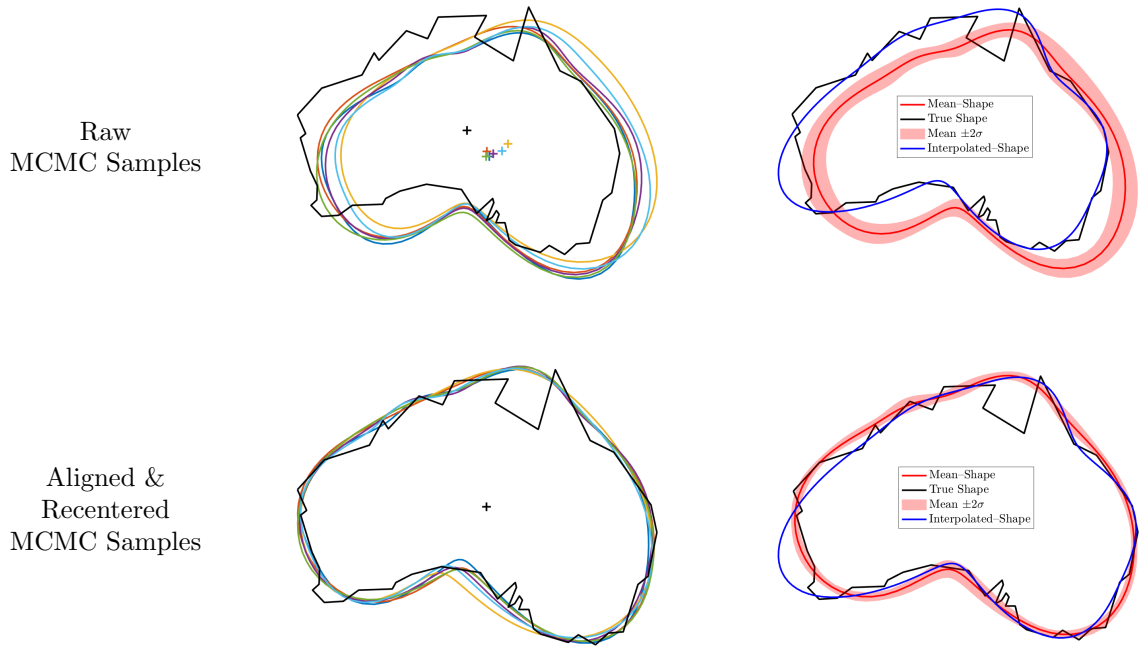
\begin{figure}
\centering
\begin{tikzpicture}
\matrix[column sep={0pt, -10pt, 10pt, 0pt}, row sep=0pt] {

    &\node{\parbox{3cm}{\centering Raw\\MCMC Samples}}; &
    &\node{\includegraphics[width=0.4\textwidth]{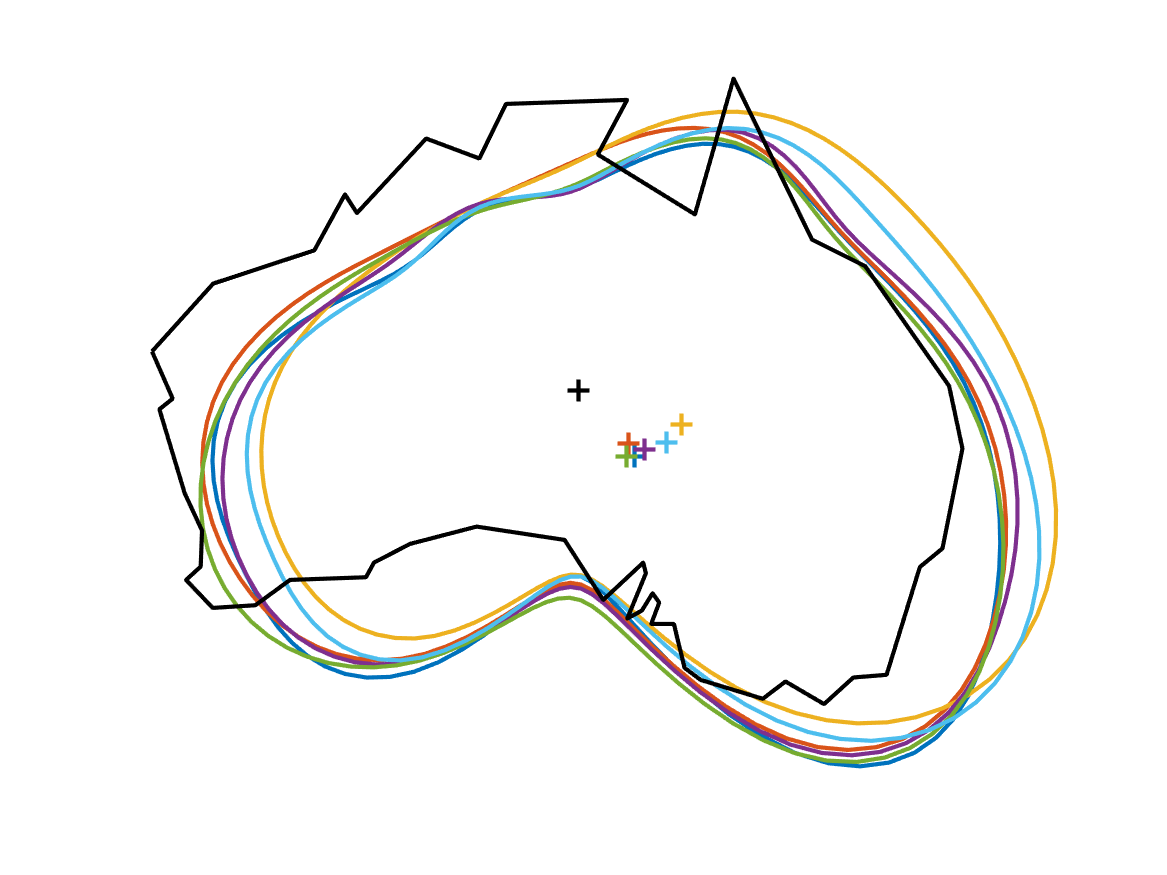}}; &
   \node{\includegraphics[width=0.4\textwidth]{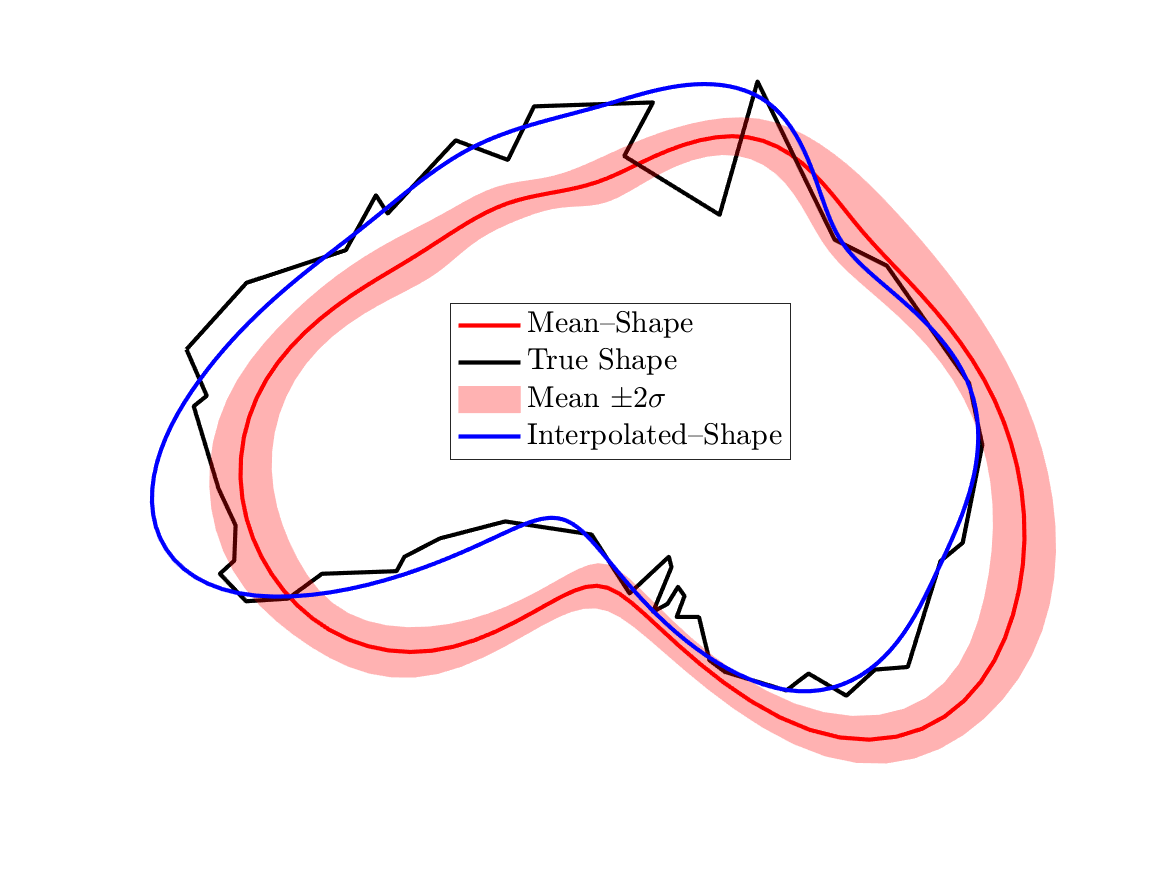}}; & \\
   &\node{\parbox{3cm}{\centering Aligned \& \\Recentered \\MCMC Samples}}; &
    &\node{\includegraphics[width=0.4\textwidth]{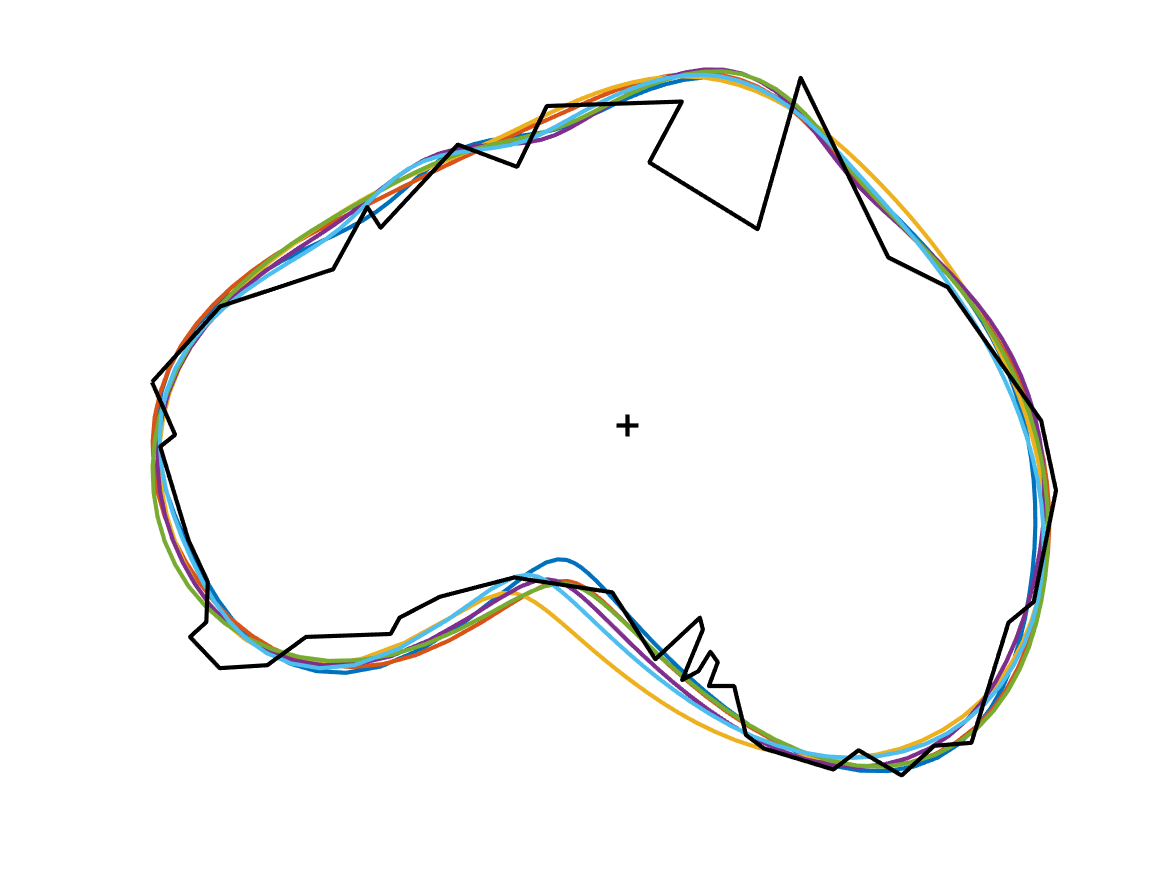}}; &
   \node{\includegraphics[width=0.4\textwidth]{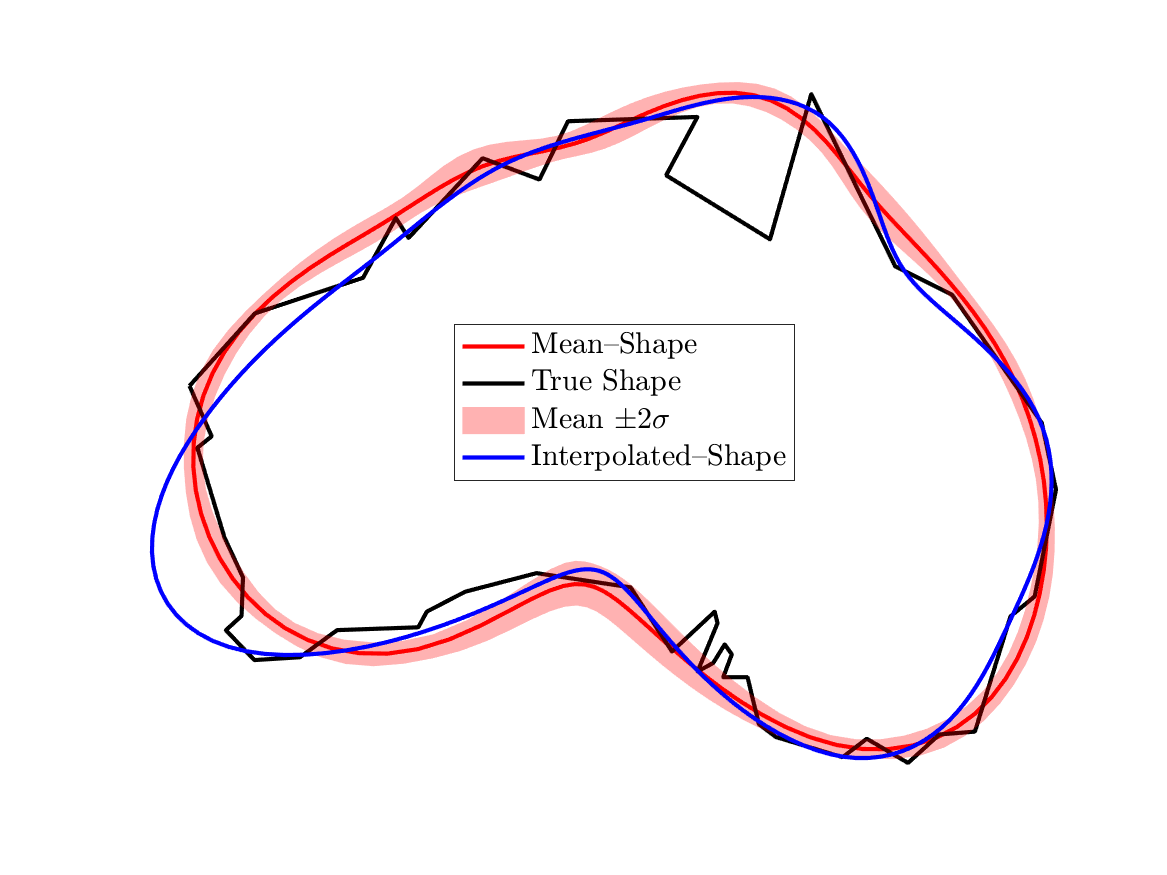}}; & \\
   };
\end{tikzpicture}
\caption{Effect of the post-processing alignment step in shape reconstruction. (Top left) 6 random reconstruction samples from the posterior distribution and their associated center of mass (plus markers), (top right) corresponding statistical inference for all samples. (Bottom left) The same 6 samples after alignment using (\ref{eq: recentering}) and (bottom right) corresponding confidence intervals. The parameters are given in Figure \ref{fig:Nspline=Nobs=12 reconstructions}.}
\label{fig:Recentering idea}
\end{figure}

Next we address an important consequence of using phaseless
far-field data.
Translations of a scatterer change the phase of its far field under illumination
by a fixed plane wave, but do not change the far-field intensity.
Consequently, our likelihood function assigns the same probability density
to translations of the same shape.
We note that translation of a shape changes its representation
$\boldsymbol{\xi}$ in our spline basis
so that, overlooking complexities of the spline basis,
multiple distinct parameter vectors $\boldsymbol{\xi}$ in our prior space
give the same data
(see for example, the samples visualized in Figure \ref{fig:Recentering idea}).

Since our goal is not to compute the position of the scatterer, but to
determine its shape, we implement a simple  alignment algorithm to
adjust for translations by aligning all of the samples about a common
origin.
The idea is to compute the centre of mass of each sampled shape,
and translate all of the samples so that they have a common centre of mass.
To that end,
we choose an arbitrary point $\mathbf{t} \in \mathbb{R}^2$ to serve as the origin of our aligned samples.
Then for each sampled shape,
we apply the translation
\begin{equation}
    \mathbf{x} \mapsto \mathbf{x} + (\mathbf{t} - \mathbf{c}), \label{eq: recentering}
\end{equation} where $\mathbf{c} \in \mathbb{R}^2 $ is the shape's center of mass (CoM).
The CoM is a simple computation, following \cite{bourke1988calculating}.
The recentering using (\ref{eq: recentering}) ensures that all of
the aligned shapes
have the same CoM. The recentering computation does not require knowledge of the
CoM of the reference shape and, in general, the CoM of the aligned shapes
may not coincide with the CoM of the reference shape.

In the numerical results section
we use the MCMC samples for uncertainty quantification by computing the
mean and variance of the radius of the corresponding aligned shapes.
For any given angle $\theta$, the radius $r(\theta)$ of the aligned shape is
easily computed as a post-processing step,
and spline representation is not required for the
log-radius of the aligned shapes.

\section{Numerical Results}\label{section:Numerical Results}

\begin{figure}
\centering
\begin{tikzpicture}
\matrix[column sep=-25pt, row sep=0pt] {
    &\node{\includegraphics[width=0.45\textwidth]{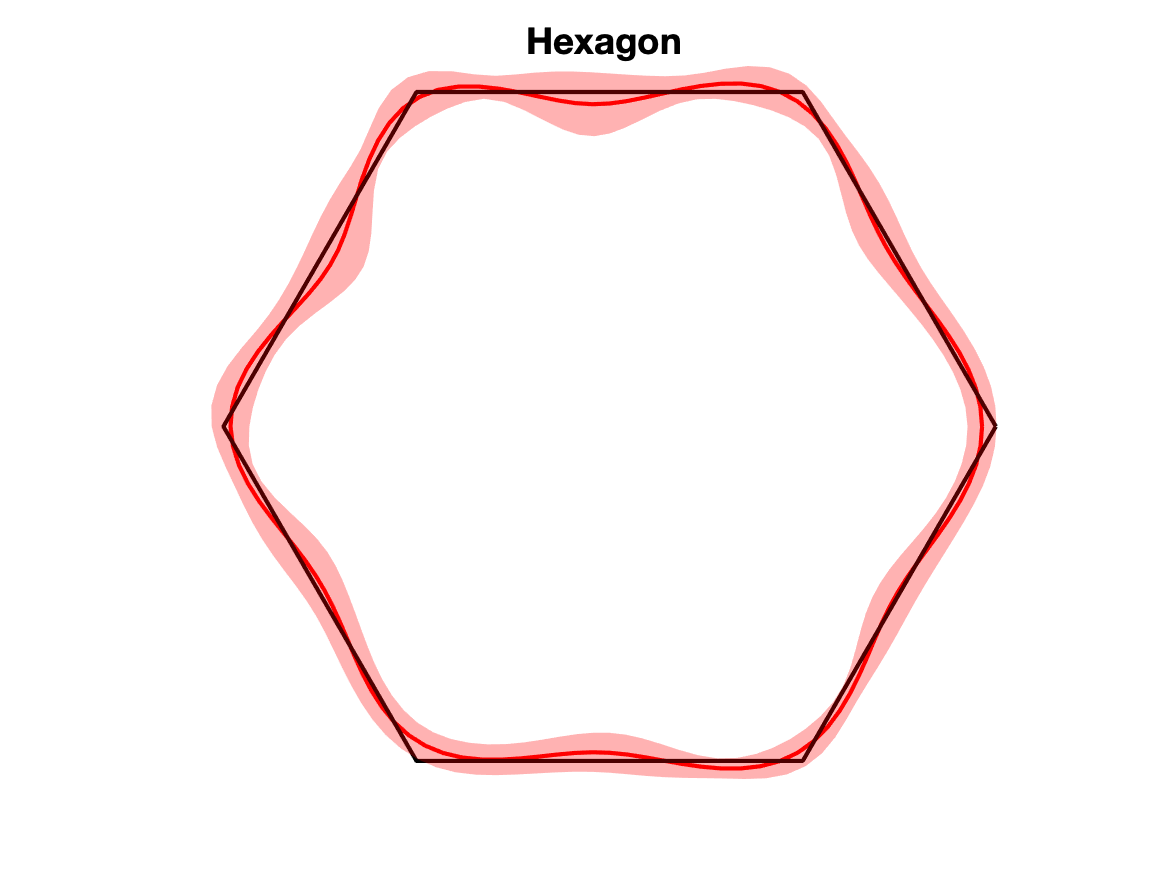}}; &
   \node{\includegraphics[width=0.45\textwidth]{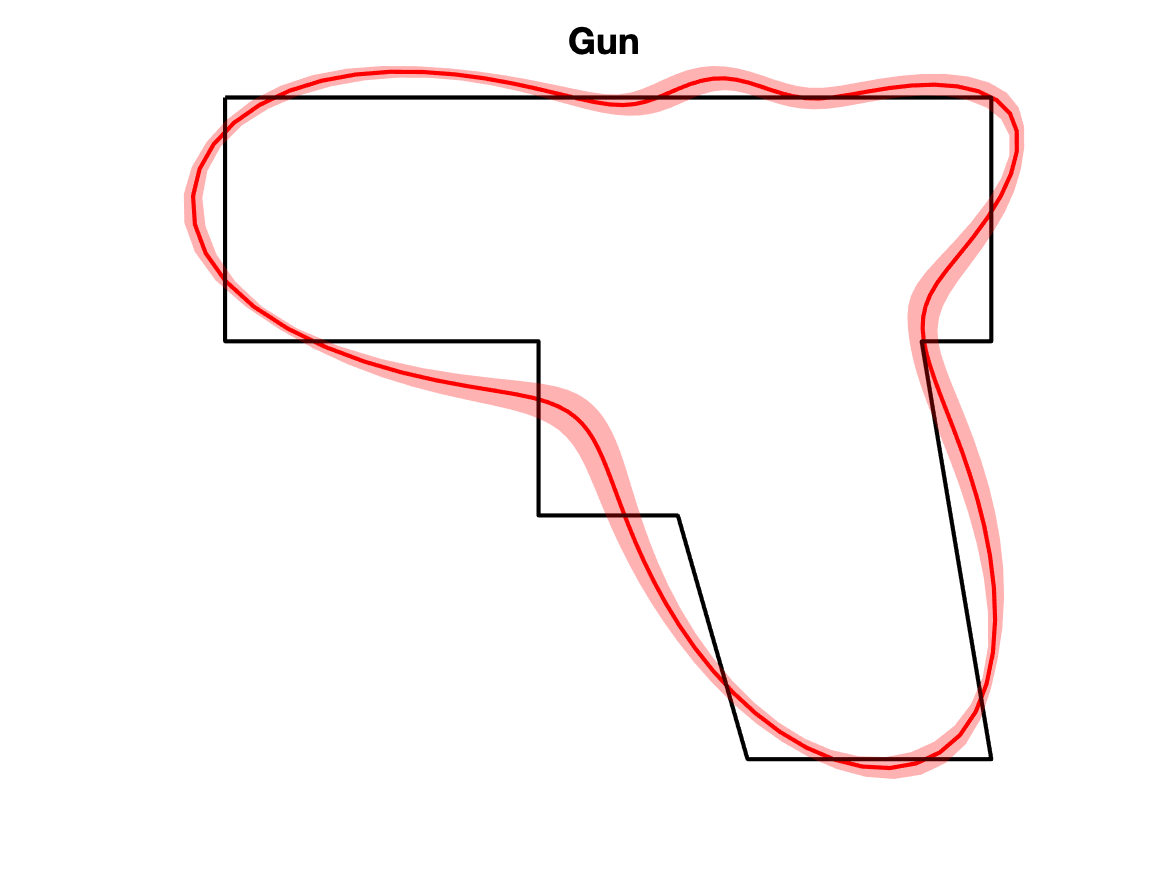}}; & \\
    &\node{\includegraphics[width=0.45\textwidth]{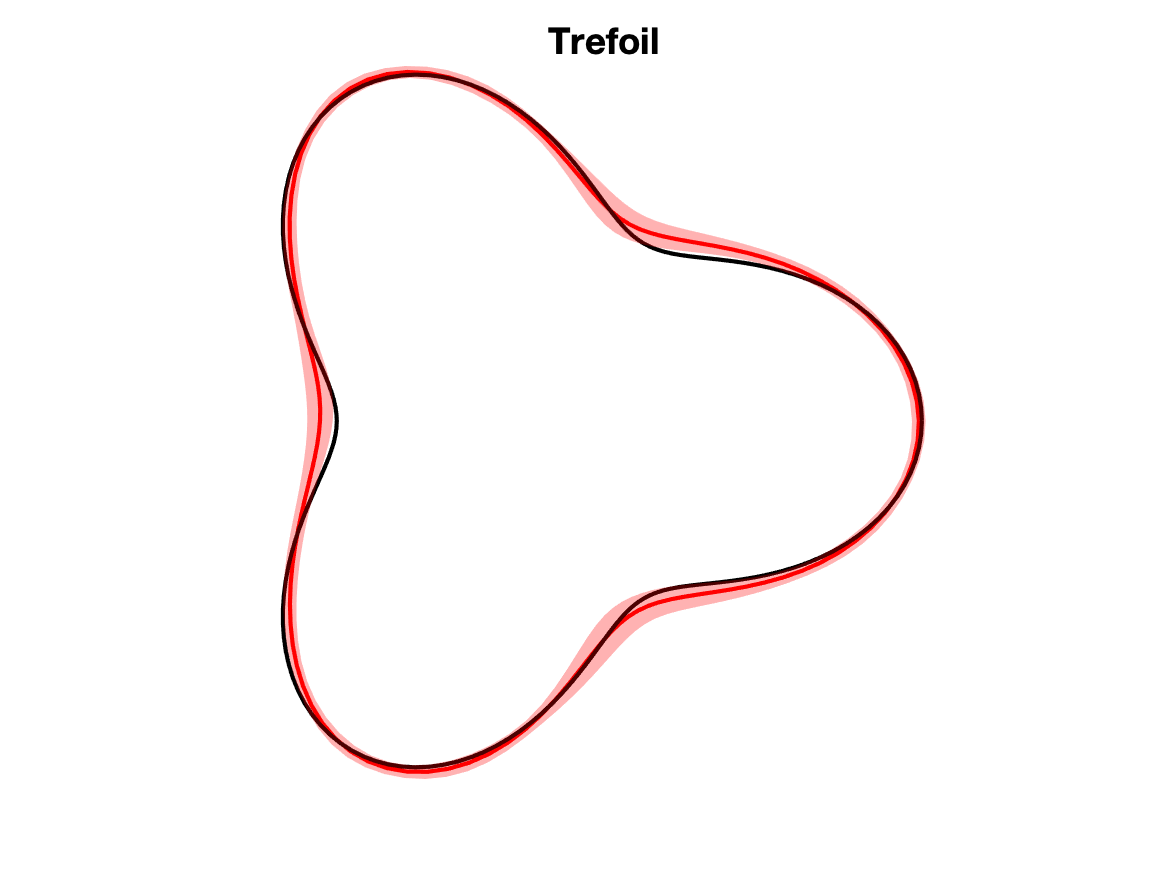}}; &
   \node{\includegraphics[width=0.45\textwidth]{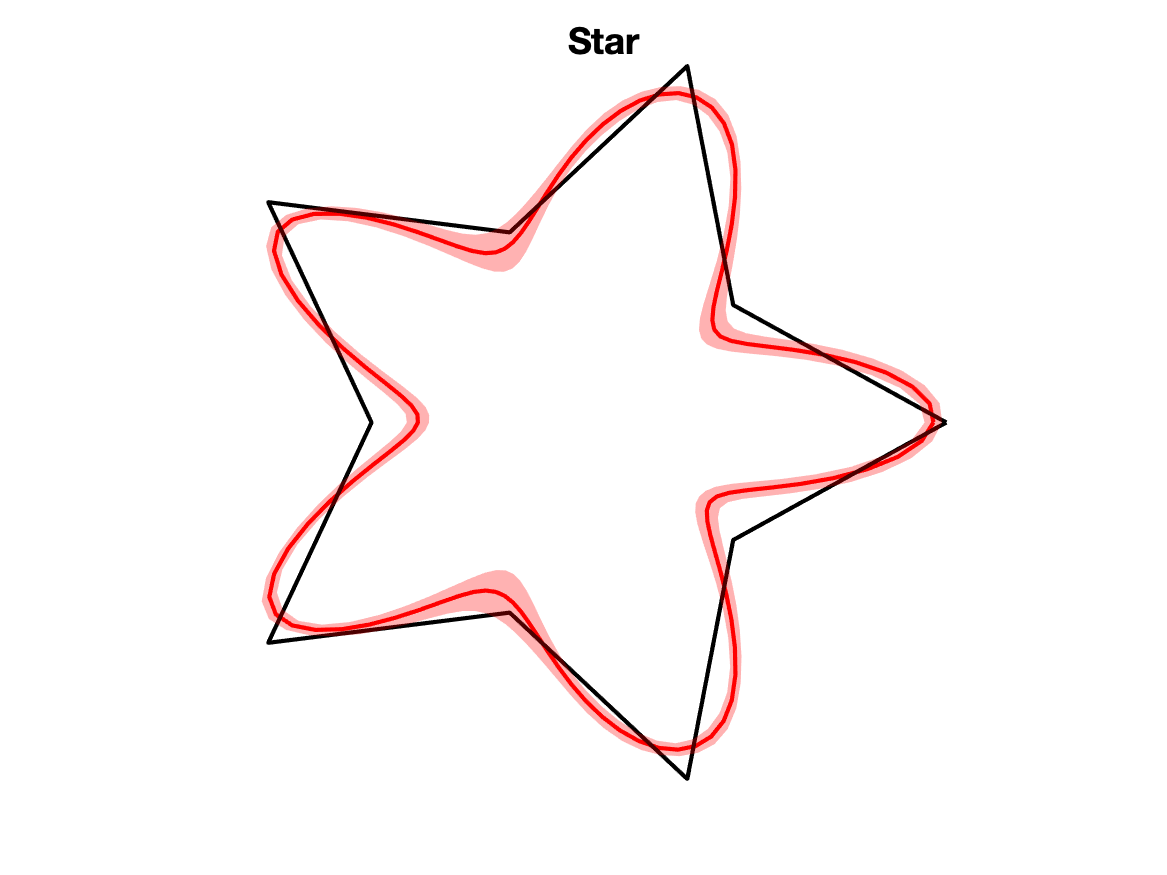}}; & \\
   &\node{\includegraphics[width=0.45\textwidth]{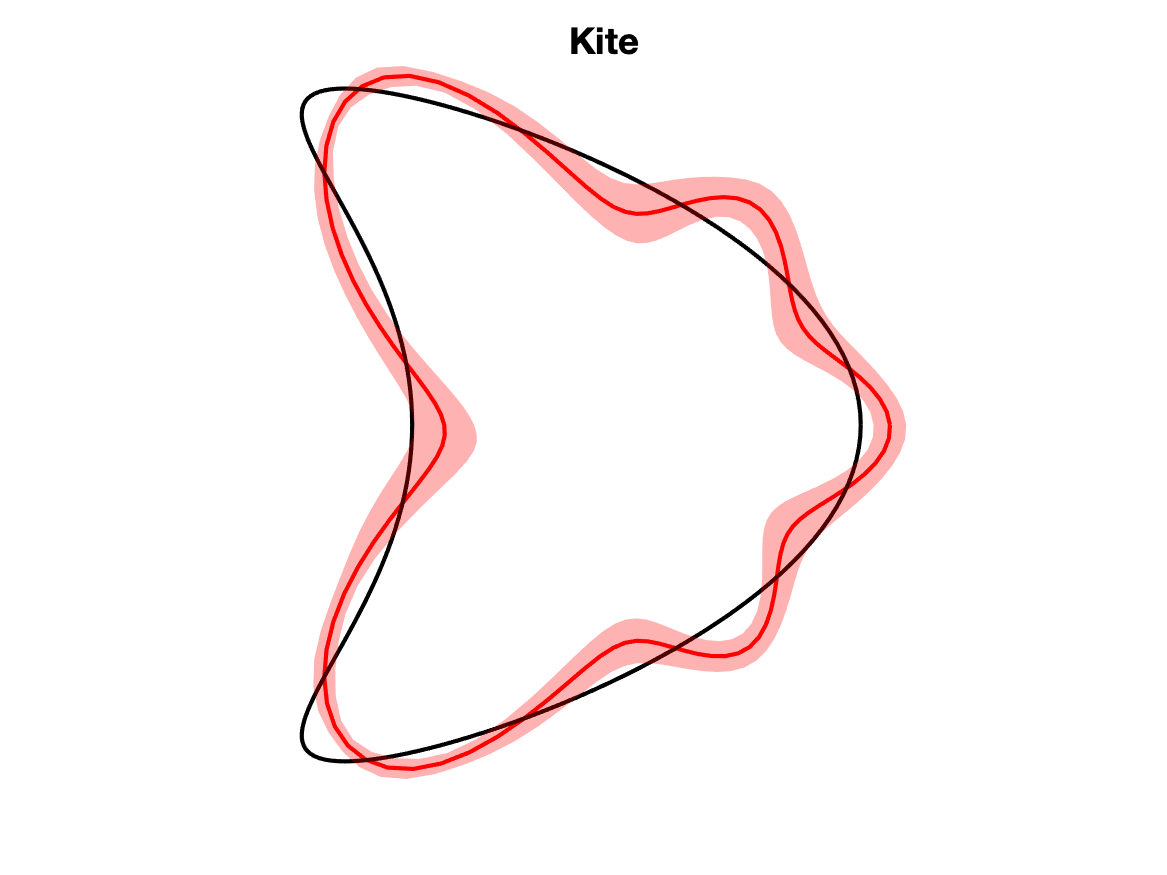}}; &
   \node{\includegraphics[width=0.45\textwidth]{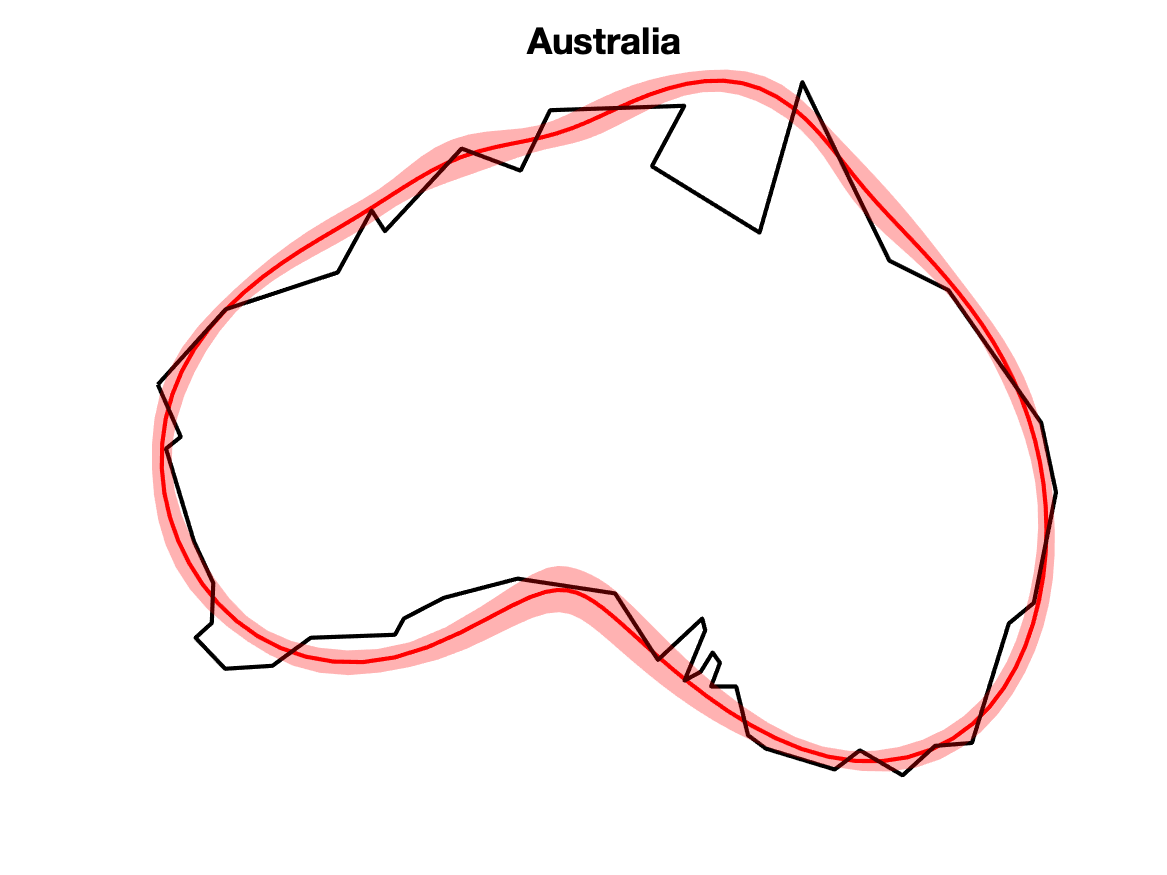}}; & \\
   };
\end{tikzpicture}
\caption{Shape reconstructions from the posterior distribution across the catalogue of shapes at $k=\pi$ with $N_\text{spline}=N_\text{inc}=N_\text{obs}=12$ and $2\%$ noise ($\hat{\sigma}=0.02$). The black lines describe the `true shapes', the red line represents the mean-shape, 
and the red-shaded area gives the associated mean $\pm 2\sigma$ confidence interval.}
\label{fig:Nspline=Nobs=12 reconstructions}
\end{figure}

We demonstrate the capability of our algorithm using simulation studies with the six test scatterers visualized in Figure \ref{fig:Nspline=Nobs=12 reconstructions}.
Our main focus in this work is to reconstruct scatterers with smooth
boundaries.
To that end, two of our test scatterers are smooth shapes, although we
emphasise that they are not constructed using splines and so they
lie outside the training distribution of our NN surrogates.
Our remaining test scatterers are polygons, selected to explore the
performance of our algorithm across varying levels of geometric complexity.

The `Trefoil' and `Hexagon' are expected to be easier targets for reconstruction, see e.g.~\cite{ganesh2020efficient} and we use them as a baseline to verify fundamental algorithm performance. The `Kite' is a standard benchmark scatterer used in \cite{colton:inverse} and further analyzed in \cite{palafox2017point, yang2020bayesian}. 
The polygonal `Gun', `Star', and `Australia' scatterers are challenging to reconstruct because they have complex shapes and discontinuous normals.
Indeed, these shapes require very large $N_\text{spline}$ values for accurate representation using interpolation. 
Furthermore, some features of the `Star' and `Gun' shapes have dimensions approaching the diffraction limit at the demonstration wavenumbers. Our numerical results demonstrate the robustness of our Bayesian algorithm even for these challenging scatterers outside our modelling assumptions.

For each test scatterer we generate our reference data by computing the far field cross section numerically, and then adding random i.i.d.~Gaussian noise. 
We use data from $N_\text{inc}$ incident waves with directions 
\begin{displaymath}
 \widehat{\mathbf{d}}_j=-(\cos \theta^\text{inc}_j, \sin \theta^\text{inc}_j),
 \quad \theta^\text{inc}_j = 2 \pi (i-1)/N_\text{inc}, \quad j=1,2,\dots,N_\text{inc},
\end{displaymath}
 and similar equispaced $N_\text{obs}$ observation directions.
 The noise $\eta_{i,j}$ is an independent sample of the Gaussian distribution with zero mean and standard deviation $\sigma=\hat{\sigma}||\mathbf{f}||_\infty$,
 where $\hat{\sigma}$ is the standard deviation of the noise
normalized by the cross section strength.
For the smooth scatterers (Trefoil and Kite) we use the Nystr\"om scheme to generate the data, but the parametrization
of the shape is wholly different to the one we use for reconstruction, so we avoid the ``inverse-crime''. For the polygonal scatterers we generate the data using
a coupled FEM-BEM scheme (see~\cite{hawkins2024efficient}), where the finite element method (FEM) is used to
accurately represent the polygonal boundary and the scattering nearby.

Our NN surrogate is a feed-forward network with a single hidden layer with $N_\text{HL}$ neurons. The hidden layer employs the hyperbolic tangent sigmoid activation function ($\tanh$) and the output layer uses linear activation. 
The prior space is as described in Section~\ref{section:inverse-shape problem}
with $-\log R_0 = \log R = 0.5$, so that, at the knots $\xi_\ell$ the  radius $r_\ell$ of the scatterers in the prior space
satisfies $r_\ell \in [0.61,1.65]$
for $\ell=1,\dots,N_\text{spline}$. 
Except otherwise stated,
in all of our experiments  $N_\text{spline}= N_\text{obs} = N_\text{inc}=12$  with the  reference wave direction
$\widehat{\mathbf{d}}_1 = (-1,0)^\top$. 
We present results for wavenumbers $k=\pi,2 \pi$ and 
different NNs are required for each wavenumber.

\begin{figure}
\centering
\begin{tikzpicture}
\matrix[column sep=0pt, row sep=0pt] {
    &\node{\includegraphics[width=0.5\textwidth]{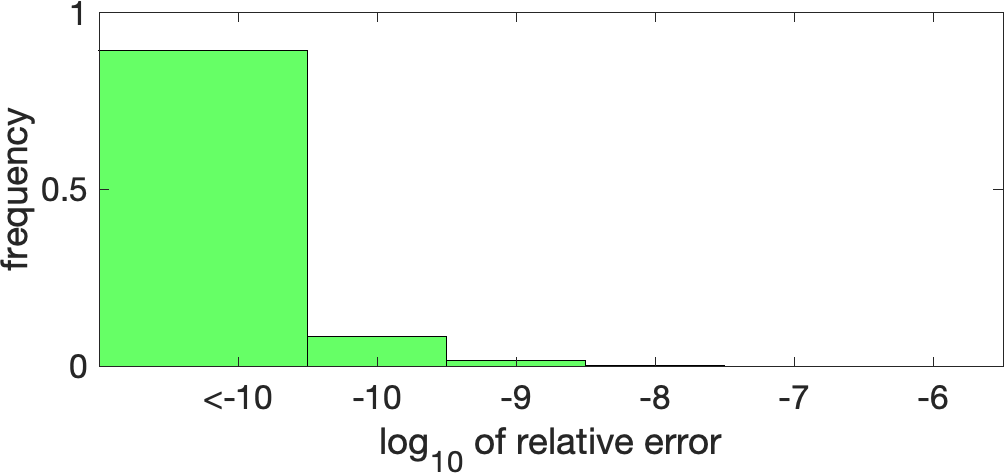}}; &
   \node{\includegraphics[width=0.5\textwidth]{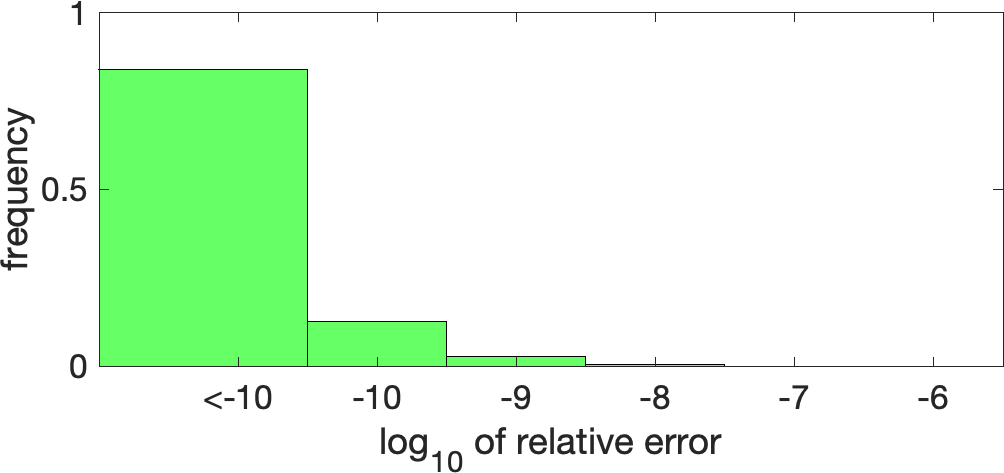}}; & \\
   };
\end{tikzpicture}
\caption{Histograms showing normalized frequency of the relative error (\ref{eq: n vs n+5 rel error}) for $1000$ random shapes $  D(\boldsymbol{\xi}^{(i)})$ with $\xi_i \sim  \mathcal{U}[-0.5,0.5]$ and $N_\text{spline}=N_\text{obs}=12$ at $k=\pi $ (left) and $k=2\pi $ (right) both with Nystr\"om discretization parameter $n=100$.}
\label{fig:Nystrom n vs n+5}
\end{figure}

Training data is generated following the Monte Carlo procedure described in Section \ref{section: forward computation}. Figure \ref{fig:Nystrom n vs n+5} shows that $\epsilon_\text{tol}=10^{-6}$ in (\ref{eq: n vs n+5 rel error}) can be achieved with Nystr\"om discretization parameter $n=100$ for both wavenumbers. The complete dataset with $M$ datapoints is therefore an input $12 \times M$ matrix and an output $12 \times M$ matrix.
For training, the data is partitioned using Matlab's default $70-15-15$ split, that is, $70\%$ for training, $15\%$ for validation during the optimization process, and $15\%$ for final testing to assess generalization performance. Both input and output matrices are preprocessed with a z-score normalization, so that the data is centered at zero mean and scales to unit variance. 
The NN with hyperbolic tangent function (tanh) activation function is trained using the Levenberg-Marquardt algorithm implemented in the Matlab Neural Network Toolbox. No explicit regularization was employed since no overfitting was observed for the chosen parameters. %

\begin{figure}
\centering
\begin{tikzpicture}
\matrix[column sep=0pt, row sep=0pt] {
    &\node{\includegraphics[width=0.5\textwidth]{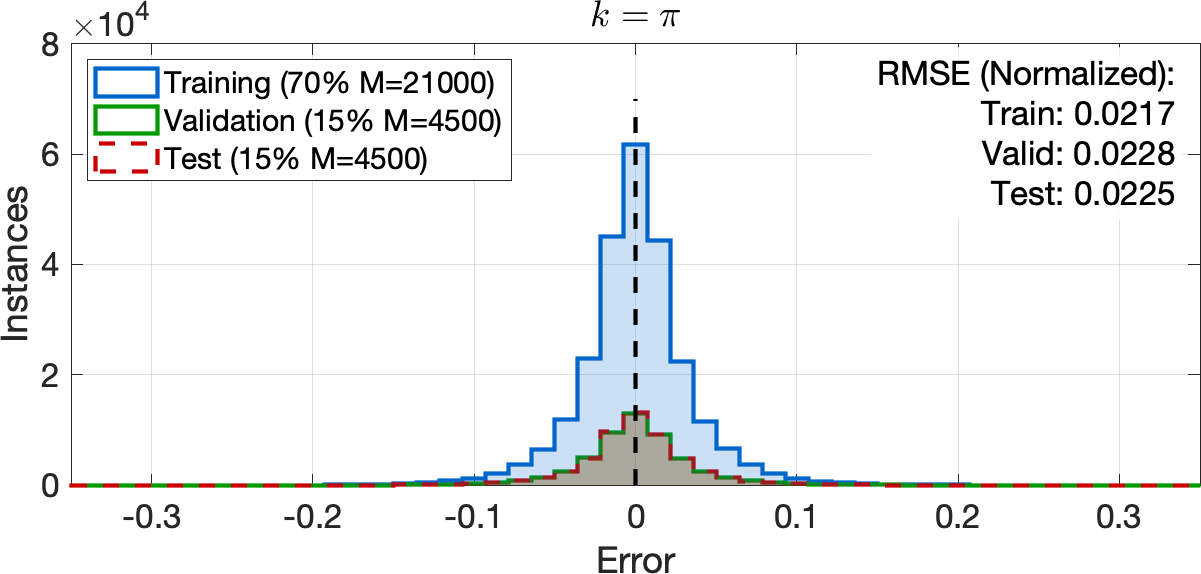}}; &
   \node{\includegraphics[width=0.5\textwidth]{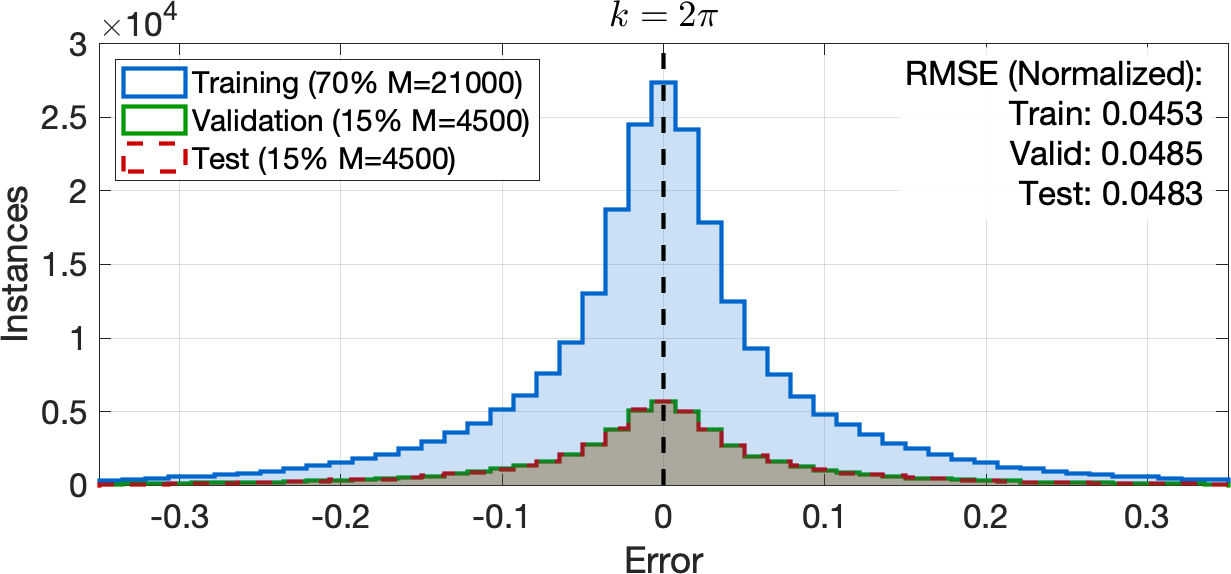}}; & \\
   };
\end{tikzpicture}
\caption{Error histogram of the NN predictions for $k=\pi $ (left) and $k=2\pi $ (right) surrogates used for reconstructions. Training (blue), 
validation (green), and test (red, dashed) sets show consistent zero-centered distributions with similar low RMSE values across all datasets.}
\label{fig:NN errors}
\end{figure}

For the surrogate NN training we choose a single hidden layer of  width $N_\text{HL}$. We use $N_\text{HL}=200,300$ for the wavenumbers $k=\pi,2\pi$ respectively.
In both cases we use $M=30,000$ training data points. Empirical testing showed that these configurations 
provide optimal approximation accuracy for their respective frequency regimes, with the higher wavenumber case requiring larger
$N_\text{HL}$ to capture the more complex scattering patterns typical at higher frequencies. The error distributions
plotted in Figure \ref{fig:NN errors} demonstrate consistent performance across training, validation, and test sets. The $k = \pi$ network achieves normalized root mean squared error (RMSE) values of 0.0217 (training), 0.0228 (validation), and 0.0225 (test), while the $k = 2\pi$ network yields RMSE values of 0.0453 (training), 0.0485 (validation), and 0.0483 (test). The zero-centered error distributions and similar RMSE values across all datasets indicate successful training without overfitting.

For all reconstruction experiments, we employ Gibbs sampling with starting point for the shape parameters $\boldsymbol{\xi}^{(0)}$ corresponding to the unit circle (with all log-radius values set to zero) within our prior space. The initial value for the noise standard deviation is set to $\sigma^{(0)} = 0.1$. We generate $15,100$ Gibbs MCMC samples for each reconstruction experiment, with a default burn-in period of $100$ samples. Since our burn-in criterion is not adaptive and the MCMC method can occasionally become trapped in parameter regions, we implement a diagnostic approach for burn-in assessment. By recording log-posterior values, we evaluate convergence through trace plot analysis. If the default burn-in proves insufficient (indicated by continued drift or slow mixing in the early portion of the chain), we discard additional samples from the beginning of the chain until satisfactory stationarity is achieved.


Because our experiments are simulation studies using synthetic data with a known
scatterer,
we are able to compute the $L_2$ error in the mean of the
aligned shapes.
In general the aligned shapes may correspond to a translation of the original
scatterer and this must be accounted for.
We do this by first
obtaining the mean-radius from the re-aligned MCMC samples (using (\ref{eq: recentering}) with $\mathbf{t}=\mathbf{0}$), then
by determining the translation vector $\mathbf{t}^\text{opt}$ that minimizes the $L_2$ distance between the reconstructed mean-shape and the true boundary: $\mathbf{t}^\text{opt} = \arg\!\min_\mathbf{t} ||\partial D - (\partial D^\text{rec} + \mathbf{t})||_2$. This optimization step is performed solely to ensure that $L_2$ error measurements reflect shape reconstruction quality and neglect translational misalignment. 

In 
Table \ref{table:cpu times} we demonstrate the substantial computational advantage achieved by our NN surrogate compared to direct Nystr\"om boundary integral evaluations by tabulating the CPU time required 
for $1000$ forward model evaluations of $\mathbf{f}^1$ on randomly sampled geometries. Using a four-core 1.4 GHz Intel i5 CPU, the Nystr\"om method with discretization parameter $n = 100$ requires approximately $290$ seconds, while the trained NN surrogate requires approximately $0.012$ seconds. This represents a $24,000\times$ speedup, making the extensive MCMC sampling required for Bayesian inference computationally feasible.

\begin{figure}
\centering
\begin{tikzpicture}
\matrix[column sep=0pt, row sep=0pt] {
    &\node{\includegraphics[width=0.5\textwidth]{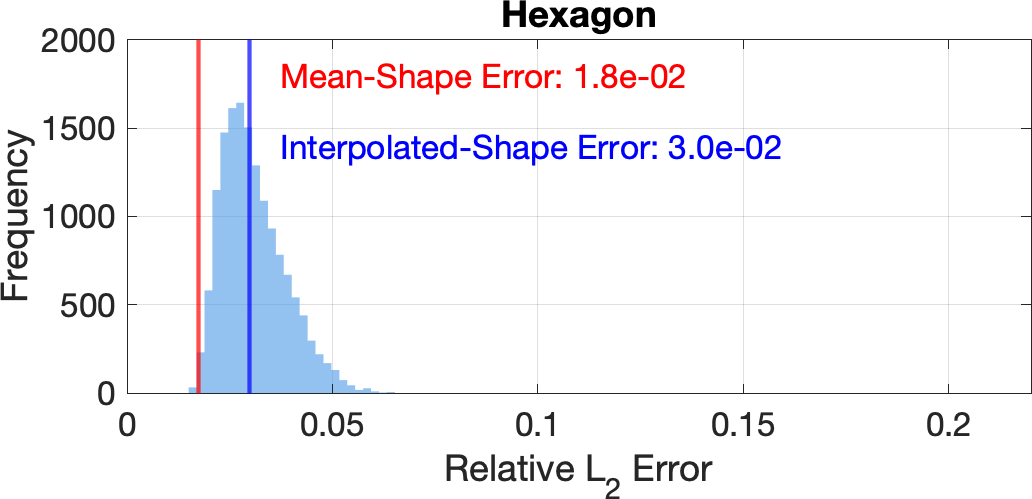}}; &
   \node{\includegraphics[width=0.5\textwidth]{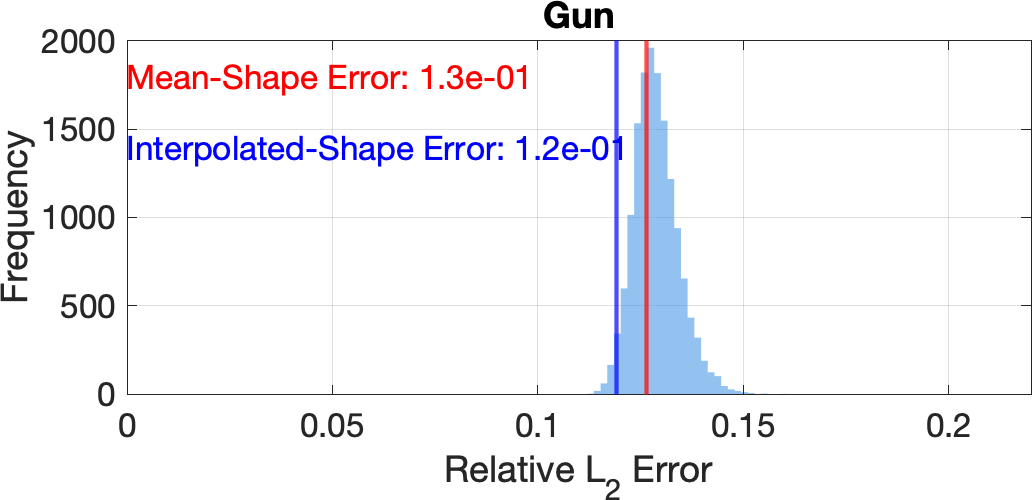}}; & \\
    &\node{\includegraphics[width=0.5\textwidth]{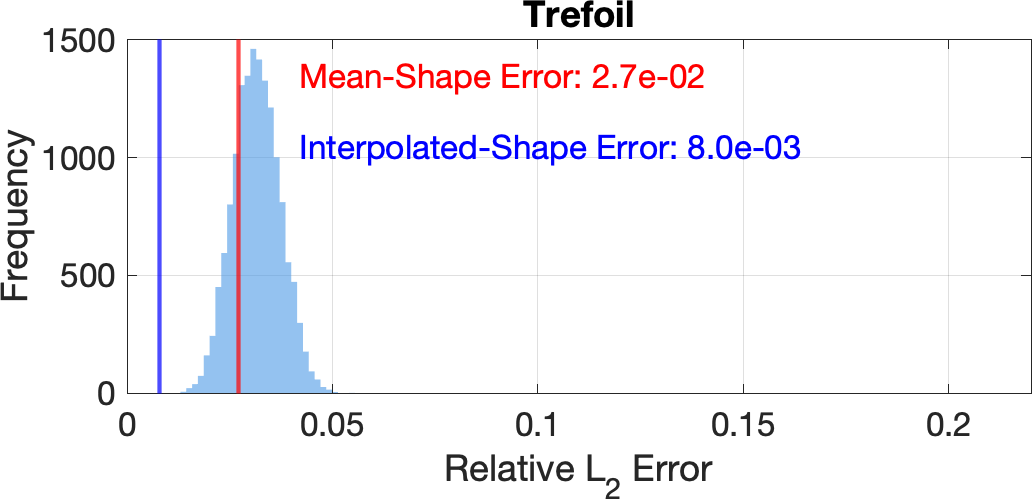}}; &
   \node{\includegraphics[width=0.5\textwidth]{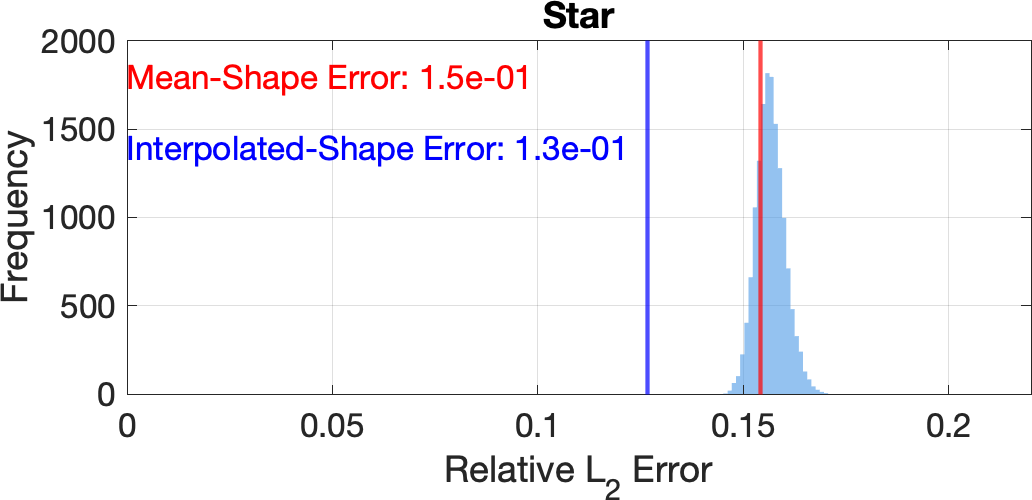}}; & \\
   &\node{\includegraphics[width=0.5\textwidth]{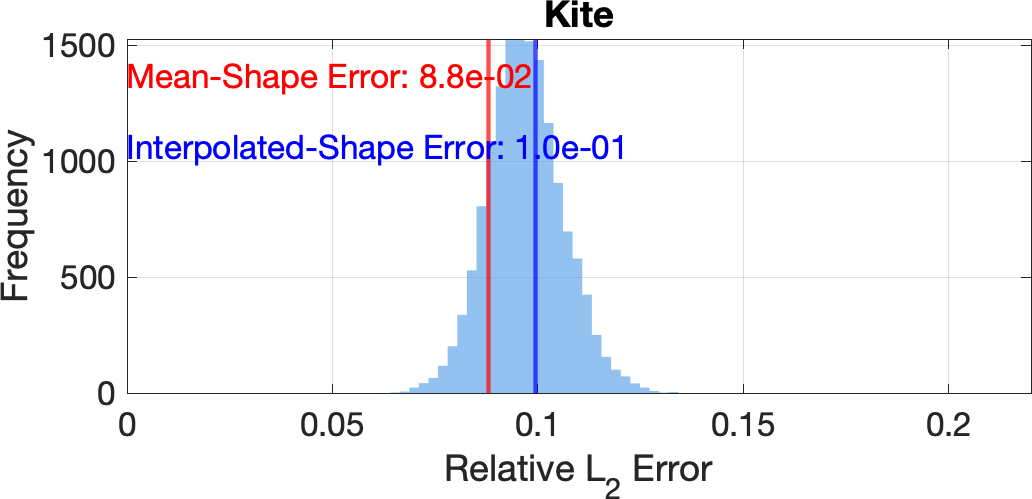}}; &
   \node{\includegraphics[width=0.5\textwidth]{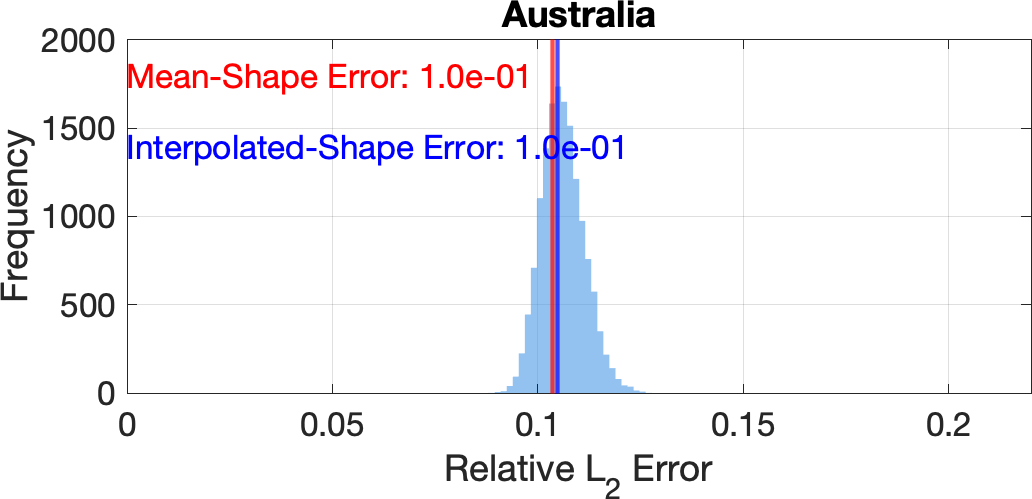}}; & \\
   };
\end{tikzpicture}
\caption{Histograms of relative $L_2$ errors for the samples in Figure \ref{fig:Nspline=Nobs=12 reconstructions} ($k=\pi$ with $N_\text{spline}=N_\text{inc}=N_\text{obs}=12$ and $2\%$ noise ($\hat{\sigma}=0.02$)). The mean-shape error is indicated in red and the interpolated-shape error is given in blue.}
    \label{fig:PER SAMPLE VS MEAN ERROR}
\end{figure}

In
Figure \ref{fig:Nspline=Nobs=12 reconstructions} we demonstrate the effectiveness of our surrogate
Bayesian algorithm across our test catalogue 
by visualizing the mean-shape  and $\pm2 \sigma$ confidence interval from all re-centred MCMC samples from our posterior distribution,
computed with wavenumber $k=\pi$  and noise $\hat{\sigma}=0.02$. 
The error distributions for the samples in Figure \ref{fig:Nspline=Nobs=12 reconstructions} is given in Figure \ref{fig:PER SAMPLE VS MEAN ERROR}. 
In the histograms we mark the error for the mean-shape and for the shape obtained by interpolating the log-radius
using our
$N_\text{spline}=12$ cubic spline basis (see Figure \ref{fig:Recentering idea}). The latter is indicative of the 
order of reconstruction accuracy that is possible with the spline basis. In a number of cases (Hexagon, Kite, Australia) the computed mean-shape outperforms the interpolated-shape in the $L_2$ norm.

\begin{figure}
\centering
\begin{tikzpicture}
\matrix[column sep=0pt, row sep=0pt] {/
    &\node{\includegraphics[width=0.45\textwidth]{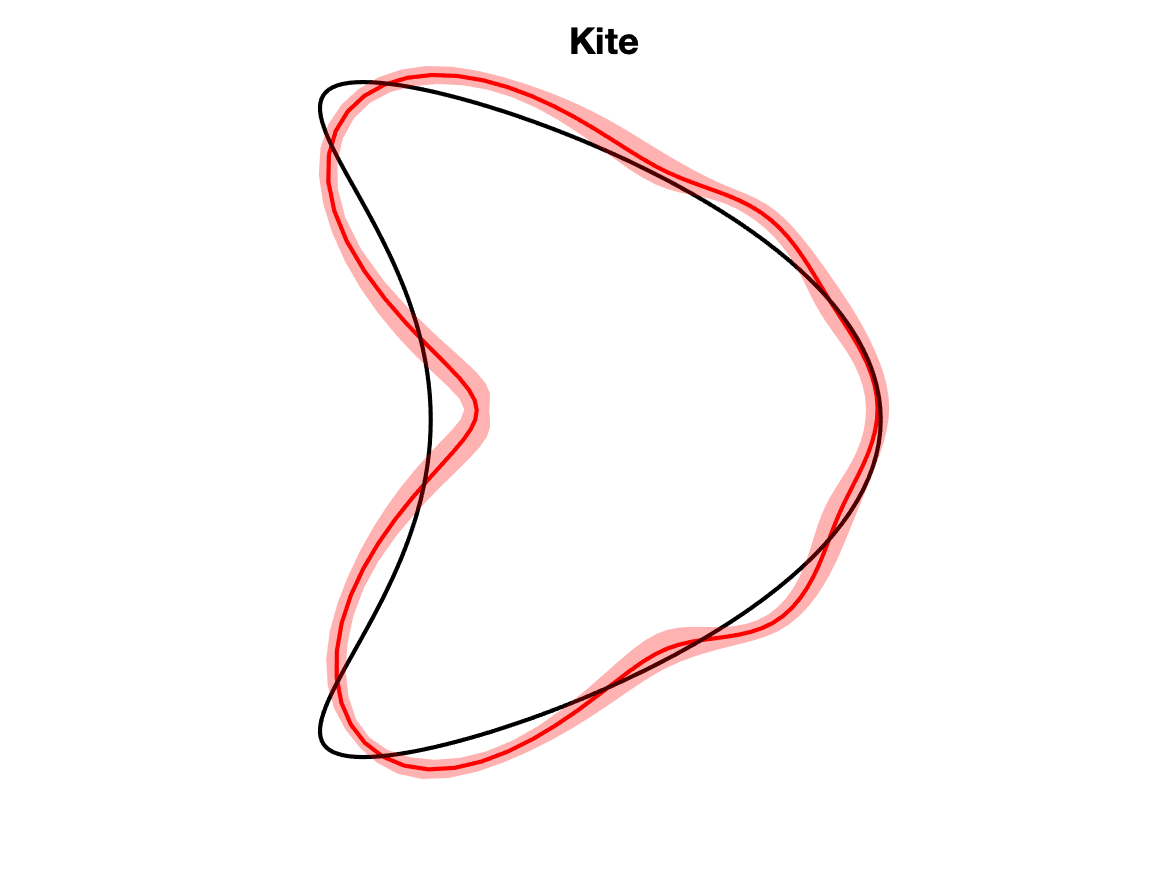}}; & 
    \node{\includegraphics[width=0.45\textwidth]{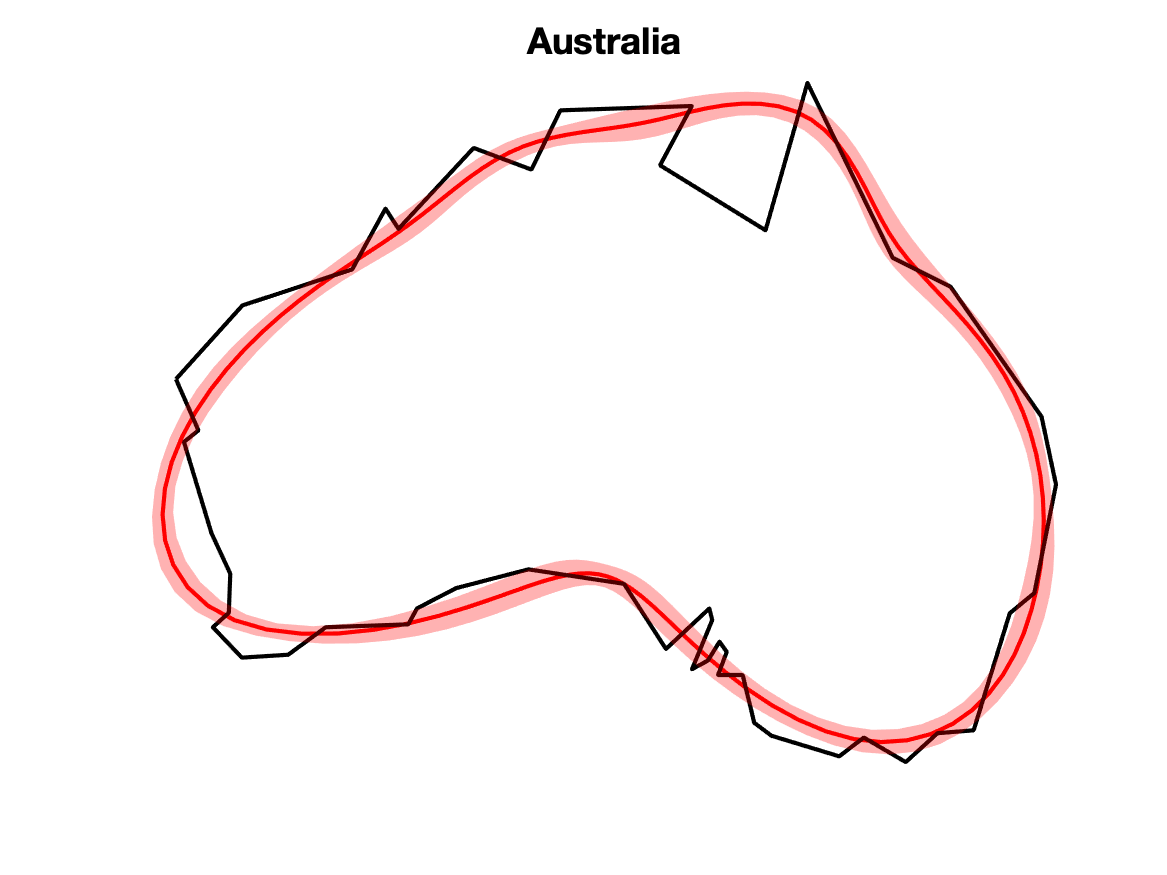}};  \\
    & \node{\includegraphics[width=0.5\textwidth]{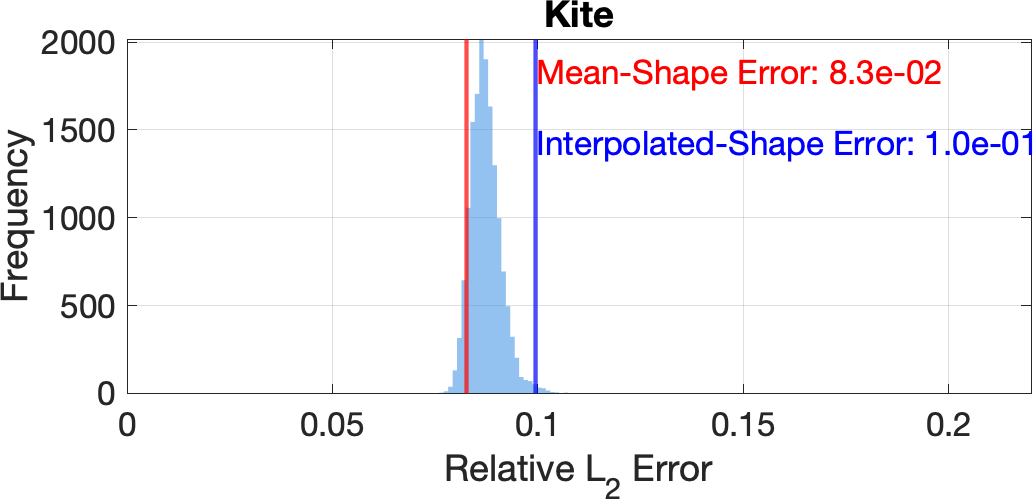}}; & 
    \node{\includegraphics[width=0.5\textwidth]{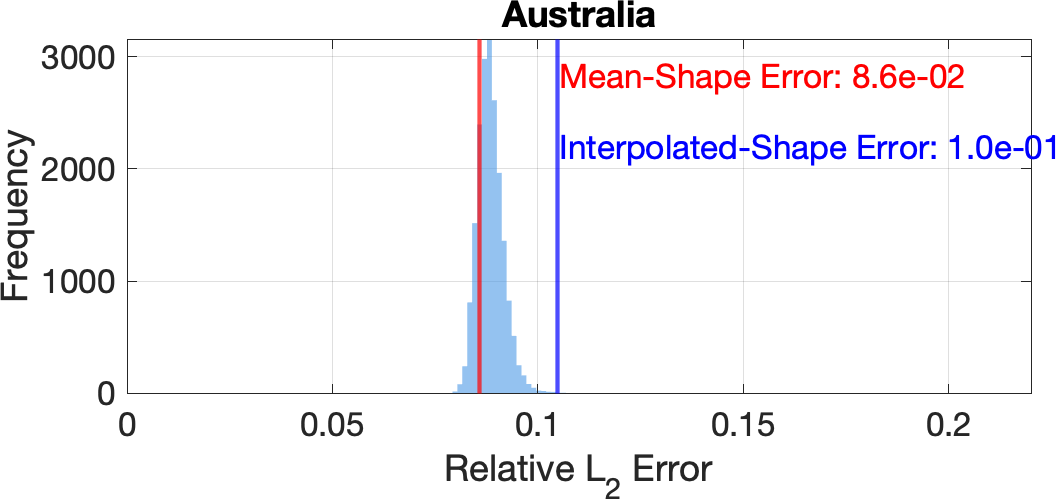}};  \\
    & \node{\includegraphics[width=0.5\textwidth]{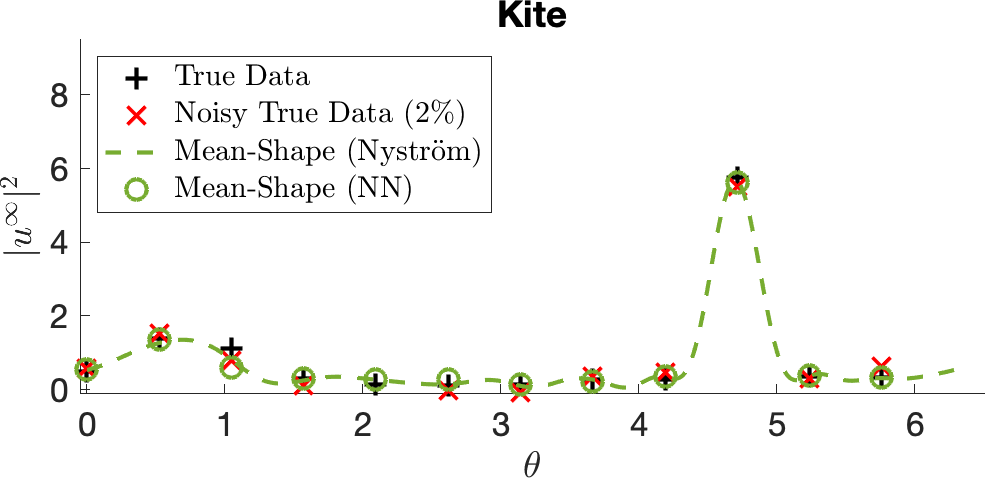}}; & 
    \node{\includegraphics[width=0.5\textwidth]{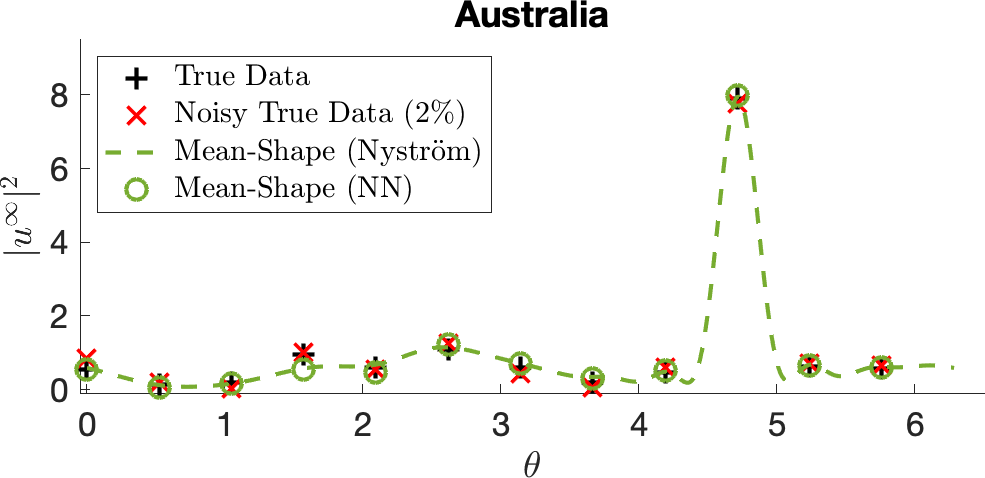}}; \\
    };
\end{tikzpicture}
\caption{Shape reconstructions (top row), $L_2$ relative error histograms (middle row) and far-field values with $\widehat{\mathbf{d}}_3=-(\cos(\pi/2),\sin(\pi/2))$ (bottom row) for the `kite' and `Australia' shapes at $k=2\pi$ with $N_\text{spline}=N_\text{inc}=N_\text{obs}=12$ and $2\%$ noise ($\hat{\sigma}=0.02$).}
\label{fig:k=2pi reconstructions}
\end{figure}

In Table \ref{table: frequency_variation} we tabulate the $L_2$ norm of the mean-shape for the scatterers in Figure \ref{fig:Nspline=Nobs=12 reconstructions}--\ref{fig:PER SAMPLE VS MEAN ERROR} at
wavenumbers $k=\pi,2\pi$. 
In Figure \ref{fig:k=2pi reconstructions} we visualize the mean-shape and $\pm2 \sigma$ confidence interval
for the Kite and Australia shapes at wavenumber $k=2\pi$ using otherwise identical parameters to
Figure~\ref{fig:Nspline=Nobs=12 reconstructions}.
Comparison between Figures~\ref{fig:Nspline=Nobs=12 reconstructions} and \ref{fig:k=2pi reconstructions}
show that the quality of the reconstructions is markedly improved at $k=2\pi$ compared with $k=\pi$, despite the norm
differences in Table~\ref{table: frequency_variation} being minor.
The qualitative differences are evident for the kite, where the oscillatory patterns near the right tip 
at $k=\pi$ are reduced
significantly at $k=2\pi$, and with narrower confidence intervals. 
This is consistent with the sensitivity analysis in \cite{neefjes2025far}. 
The error distributions for $k=2\pi$ are visualized in Figure~\ref{fig:k=2pi reconstructions}.
Comparison with Figure~\ref{fig:PER SAMPLE VS MEAN ERROR} shows that the mean error is reduced at $k=2\pi$
compared to $k=\pi$ and the standard deviation is reduced.

In the bottom row of Figure \ref{fig:k=2pi reconstructions} we visualize the data
and true far-field values for the incident wave $\widehat{\mathbf{d}}_3=-(\cos(\pi/2),\sin(\pi/2))$.
The close fit between the Nystr\"om solution (computed with $n=100$), the mean-shape and the true data demonstrates that the posterior distribution captures the true shape, and also the fidelity of the surrogate model for the mean-shape. The mean-shape boundary for the Nystr\"om computation here is represented via a high resolution cubic spline interpolant with $120$ knots.


In 
Table \ref{table:Ninc variation} 
we demonstrate the effect of the number of incident directions 
by tabulating the $L_2$ error of the mean-shape for between one and 12 incident directions for
wavenumber $k=\pi$ with and noise $\hat{\sigma}=0.02$.
In general, the results improve with increasing numbers of incident waves.
For some of the more complex shapes (Gun and Star) the limited improvement indicates that their reconstruction is 
limited by geometric complexity rather than data sufficiency.

In Figure \ref{fig:Ninc, Noise effect} we examine the effect of the number of incident waves in more detail with a histogram
of $r(\theta)$ for the Australia shape at $\theta=2.93$ computed from our MCMC samples.
Qualitatively similar results were obtained at other observation angles and for the other scatterers.
The results demonstrate concentration of the posterior distribution as the number of incident directions increases.
This plot demonstrates why the mean-shape alone cannot fully characterise the posterior. 
In particular, it is possible for the mean-shape error to be smaller than the error of each of the samples.

In Table \ref{table:Noise_variation} we demonstrate the effect of measurement  noise by tabulating the $L_2$ error
of the mean-shape for measurement noise up to 8\% with $N_\text{inc}=12$ and wavenumber $k=\pi$. 
In Figure \ref{fig:Ninc, Noise effect} (right) we examine this effect in more detail 
with a histogram
of $r(\theta)$ for the Australia shape at $\theta=2.93$ computed from our MCMC samples (similar to the left panel, which
is explained above).
We include results with $0\%$ noise as a baseline that is limited only by physical constraints (including the modeling assumptions) and the forward surrogate. 
As in the case of the number of incident waves, the limited improvement in the $L_2$ error as the noise decreases
for the more complex shapes (Gun and Star)
indicates that their reconstruction is limited by geometric complexity rather than limitations of the data.
This is expected because these are challenging shapes with non-smooth boundaries that are hard to represent using splines.

The Kite shape results in Table \ref{table:Noise_variation} warrant further consideration.
For this scatterer the mean-shape reconstruction error unexpectedly decreases with noise. 
This is a case where the mean-shape error does not capture the complexity of the posterior. In particular,
in this case the capability of the $N_\text{spline} = 12$ basis to represent the scatterer is limited, so that
in the $0\%$ noise case the posterior distribution is highly concentrated but the mean-shape is not correct (see Figure \ref{fig:Nspline=Nobs=12 reconstructions}).
In contrast, in the $8\%$ noise case the posterior is more spread, allowing the mean-shape to align better with the true boundary. 

In Table \ref{table: Nobs_variation} we demonstrate the effect of the number of observation directions by tabulating the 
$L_2$ error of the mean-shape for $N_\text{obs}$ between 3 and 12. Here we use $N_\text{inc}=3$ incident wave directions,
wavenumber $k=\pi$ and noise $\hat{\sigma}=0.02$. 
These results demonstrate the necessity of having at least $N_\text{obs}=6$ observation direction for the parameters considered.

\begin{figure}
\centering
\begin{tikzpicture}
\matrix[column sep=0pt, row sep=0pt] {
    &\node{\includegraphics[width=0.5\textwidth]{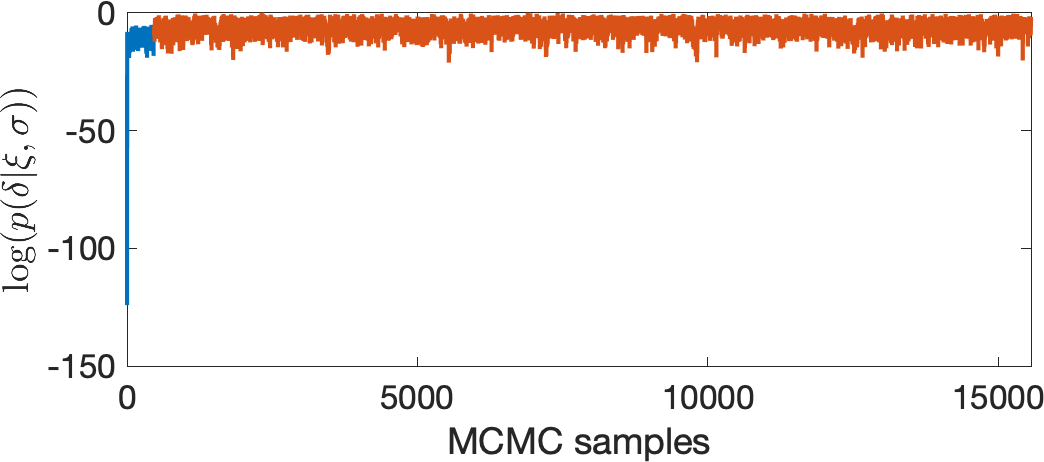}}; &
   \node{\includegraphics[width=0.5\textwidth]{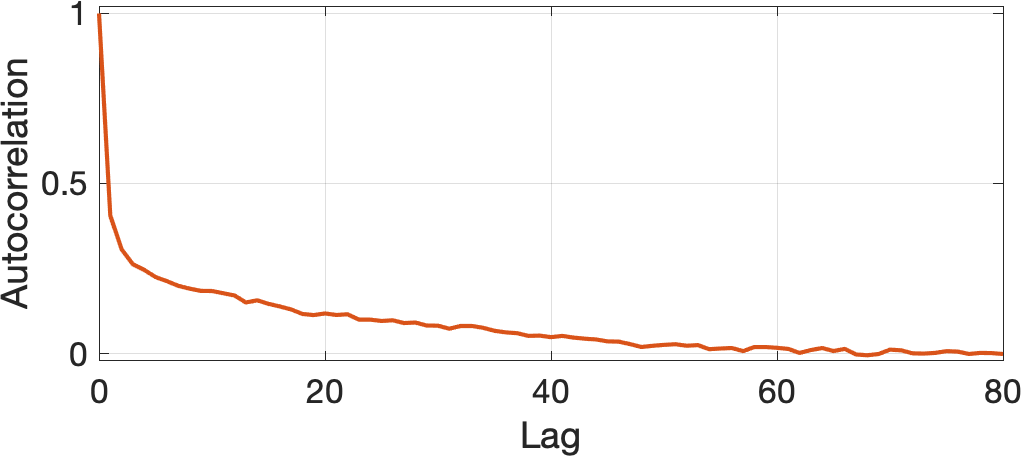}}; & \\
   };
\end{tikzpicture}
\caption{MCMC sampling diagnostics for Australia shape with $N_\text{spline}=N_\text{obs}=N_\text{inc}=12$, $k=\pi$ and $2\%$ noise ($\hat{\sigma}=0.02$). (Left) Trace plot of the log-posterior distribution across Markov Chain samples, (right) autocorrelation of the log-posterior vs lag.
}
\label{fig:autocorrelation}
\end{figure}

Finally, in Figure \ref{fig:autocorrelation} (left) we demonstrate that the Markov Chain has reached a stationary distribution
by plotting the trace of the log-posterior for the Australia shape with wavenumber $k=\pi$ and $\hat{\sigma}=0.02$.
Here the chosen burn-in of 500 is indicated by the blue.
Good mixing is indicated by the rapid fluctuations around the mean value. 
In Figure \ref{fig:autocorrelation} (right) we demonstrate the independence of the Markov Chain samples using an  
Autocorrelation analysis. 
Efficient decorrelation of the samples is indicated by rapid decay of the autocorrelation exhibits, dropping to approximately $0.5$ within $5$-$10$ lags and approaching zero by lag $20$-$30$.
The effective sample size in this case is $1270$. Qualitatively similar behavior was observed for all other experiments.
This analysis was performed using \texttt{UWerr.m} code \cite{wolff2004monte}.

\begin{table}[]
    \centering
    \begin{tabular}{c||c}
    \hline
    & CPU time \\ \hline
Nystr\"om method $n=100$ & $290 \sec$ \\
NN surrogate  &  $0.012 \sec$  \\
\hline
    \end{tabular}
    \caption{Approximate CPU times for $1000$ evaluations of the forward map $\mathbf{f}^1$ for random shapes with boundaries $ \partial D(\boldsymbol{\xi}^{(i)})$ following $\xi_i \in \mathcal{U}[-0.5,0.5]$ and $N_\text{spline}=N_\text{obs}=12$. A $24,000 \times$ speed-up is obtained using the NN surrogate.}
    \label{table:cpu times}
\end{table}

 \begin{table}[]
     \centering
     \begin{tabular}{c||cc}
     \hline
     & \multicolumn{2}{c}{$L_2$ relative error} \\
     & $k=\pi$ & $k=2\pi$ \\ \hline
 Hexagon  &  0.02  &   0.02 \\
 Gun  &  0.13  &   0.14 \\
 Trefoil  &  0.03  &   0.03 \\
 Star  &  0.15  &   0.15 \\
 Kite  &  0.09  &  0.08 \\
 Australia  &  0.10  &  0.09 \\
 \hline
     \end{tabular}
     \caption{Mean-shape reconstruction $L_2$ relative error with $N_\text{spline}=N_\text{obs}=N_\text{inc}=12$ and $2\%$ noise ($\hat{\sigma}=0.02$) for the two wavenumbers $k$ considered in our surrogates.}
     \label{table: frequency_variation}
\end{table}

\begin{table}[]
\begin{center}
    \begin{tabular}{c||cccc}
    \hline
    & \multicolumn{4}{c}{$L_2$ relative error} \\
      &  $N_\text{inc}=1$    &  $N_\text{inc}=2$    &  $N_\text{inc}=4$    &   $N_\text{inc}=12$ \\ \hline 
Hexagon &  0.10   &  0.06   &  0.03   &  0.02 \\
Gun  &  0.17  &  0.16  &   0.15   &  0.13 \\
Trefoil  &  0.18   &  0.06  &  0.03  &  0.03 \\
Star  &  0.16  &  0.15   &  0.15  &  0.15 \\
Kite  &  0.27  &  0.14   &  0.09   &  0.09 \\
Australia  &  0.15  &  0.12  &   0.11   &  0.10 \\
\hline
    \end{tabular}
    \caption{Mean-shape reconstruction $L_2$ relative error with $N_\text{spline}=12=N_\text{obs}$, $k=\pi$ and $2\%$ noise ($\hat{\sigma}=0.02$) for increasing number of incident waves. In all cases we ensure that $\widehat{\mathbf{d}}_{1} = 0$.}
\label{table:Ninc variation}
\end{center}
\end{table}

\begin{table}[]
    \centering
    \begin{tabular}{c||cccc}
    \hline
    & \multicolumn{4}{c}{$L_2$ relative error}\\
   & Noise 8\% & Noise 4\% & Noise 2\% & Noise 0\% \\
  \hline 
Hexagon  &  0.05    &  0.03    &  0.02    &  0.02 \\
Gun  &  0.14    &  0.13    &  0.13    &  0.13 \\
Trefoil  &  0.06    &  0.04    &  0.03    &  0.01 \\
Star  &  0.15    &  0.15    &  0.15    &  0.16 \\
Kite  &  0.07    &  0.07    &  0.09   &  0.10 \\
Australia  &  0.14    &  0.12    &  0.10    &  0.11 \\
\hline
    \end{tabular}
    \caption{Mean-shape reconstruction $L_2$ relative error with $N_\text{spline}=N_\text{obs}=N_\text{inc}=12$, $k=\pi$ for decreasing levels of noise $\hat{\sigma}$ in the available data.}
    \label{table:Noise_variation}
\end{table}

\begin{table}[]
    \centering
    \begin{tabular}{c||ccc}
    \hline
    & \multicolumn{3}{c}{$L_2$ relative error}\\
           & $N_\text{obs}=3$ & $N_\text{obs}=6$ & $N_\text{obs}=12$ \\ \hline
Hexagon  &  0.10   &  0.05  &   0.06 \\
Gun  &  0.25   &  0.13  &   0.13 \\
Trefoil  &  0.06   &  0.04  &   0.04 \\
Star  &  0.25   &  0.15  &   0.15 \\
Kite  &  0.14   &  0.12  &   0.09 \\
Australia  &  0.30   &  0.14  &   0.13 \\
\hline
    \end{tabular}
    \caption{Mean-shape reconstruction $L_2$ relative error with $N_\text{spline}=12$, $N_\text{inc}=3$ (with $\theta^\text{inc}_i=2\pi/3,4\pi/3,2\pi$), $k=\pi$ and $2\%$ noise ($\hat{\sigma}=0.02$) for increasing number of observation directions $N_\text{obs}$. }
    \label{table: Nobs_variation}
\end{table}

\begin{figure}
\centering
\begin{tikzpicture}
\matrix[column sep=0pt, row sep=0pt] {
    &\node{\includegraphics[width=0.5\textwidth]{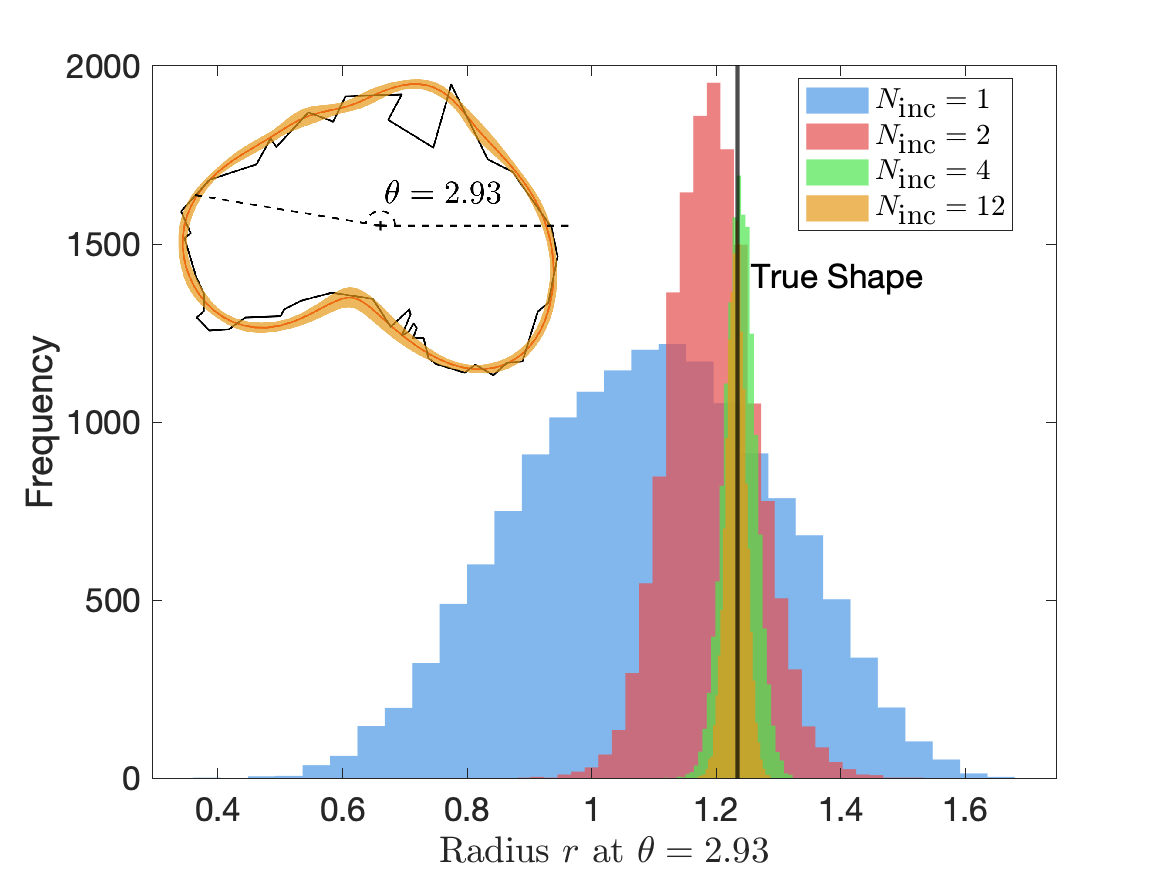}}; &
   \node{\includegraphics[width=0.5\textwidth]{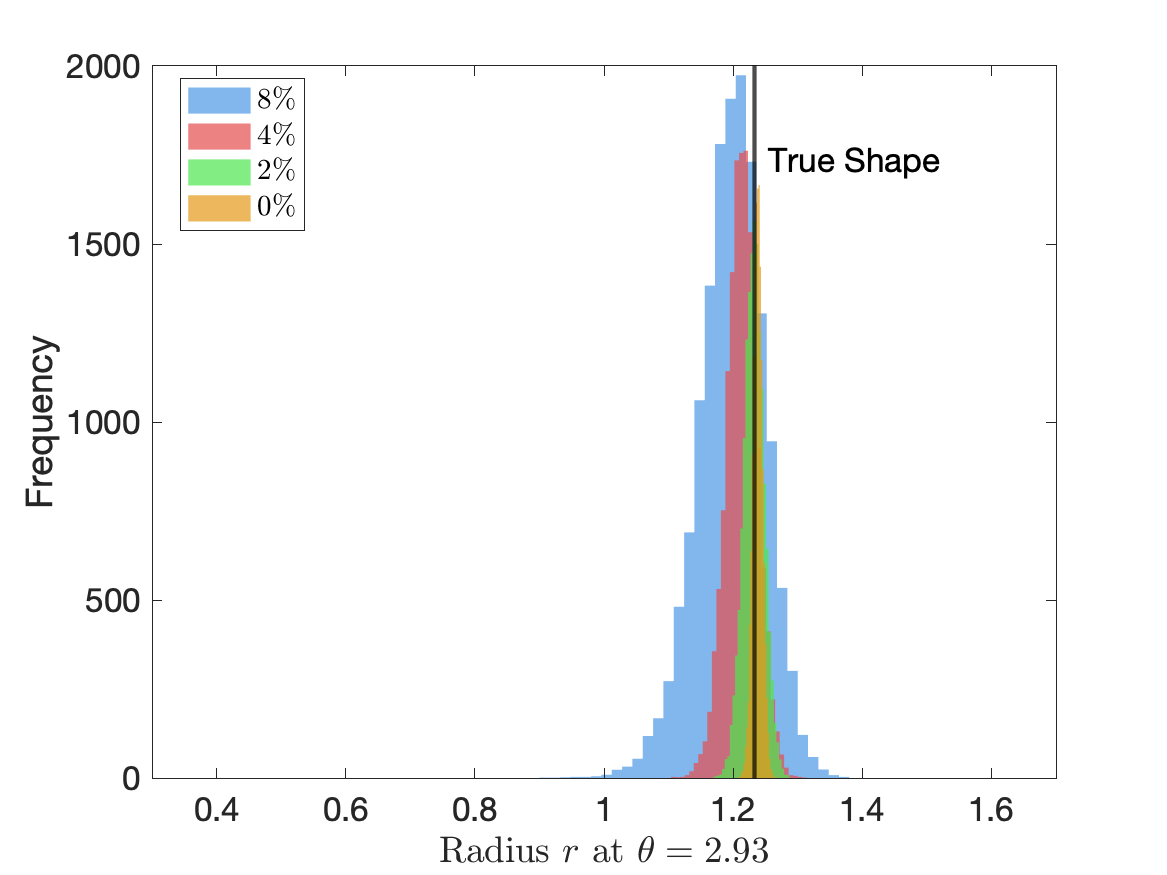}}; & \\
   };
\end{tikzpicture}
\caption{Histograms for the radii of the reconstructed samples from the posterior distribution at fixed angle  $\theta=2.93$ for Australia shape as we increase $N_\text{inc}$ (left) and noise level (right). We fix $N_\text{spline}=12=N_\text{obs}$, $k=\pi$ and $\hat{
\sigma}=0.02$ (left) and $N_\text{inc}=12$ (right). Note that the bin size of the various histograms is not constant.}
\label{fig:Ninc, Noise effect}
\end{figure}

\section*{Acknowledgements}
The authors gratefully acknowledge the support of
the Australian Research Council (ARC) through a Discovery Project Grant (DP220102243).
The work was carried out during EGN's employment as  a postdoctoral fellow under the
ARC grant awarded to SCH and MG.
In addition, MG is supported by the Simons Foundation through Grant No.~518882.

\section*{Data availability}

Far field data used for reconstructions will be made publicly available after publication.


\section*{References}
\begingroup
\renewcommand{\section}[2]{}

\bibliographystyle{unsrt}
\bibliography{Bibliography}

\end{document}